\documentclass{amsart}

\usepackage{amsmath}%
\usepackage{amsfonts}%
\usepackage{amssymb}%
\usepackage{graphicx}
\usepackage{epstopdf}
\usepackage{placeins} 
\usepackage{mathtools}
\mathtoolsset{showonlyrefs}
\usepackage{subfigure}
\usepackage{color}
\usepackage{dsfont}
\usepackage{bbm}
\usepackage{enumerate}
\usepackage{mathrsfs}

\def\calL{\mathcal{L}}

\def\calP{\mathcal{P}}
\def\calN{\mathcal{N}}

\newcommand{\Th}{{\mathcal {T}_h}}

\def\calLs{\calL_{\rm stab}}

\def\R{{\mathbb R}}
\def\N{{\mathbb N}}
\def\bV{{\mathbb V}}

\newcommand{\tHs}{{\widetilde H}^s(\Omega)}

\newcommand{\eps}{\varepsilon}
\newcommand{\vPhi}{{\mathbf \Phi}}
\newcommand{\vPsi}{{\mathbf \Psi}}
\renewcommand{\div}{{\rm div}} 
\newcommand{\grad}{{\bf grad}} 
\newcommand{\divs}{{\rm div}^s} 
\newcommand{\grads}{{\bf grad}^s}

\newcommand{\iinte}{{\mathbf \Pi}_h}
\newcommand{\inte}{{\Pi}_h}

\newcommand{\vertiii}[1]{{\left\vert\kern-0.25ex\left\vert\kern-0.25ex\left\vert #1 
    \right\vert\kern-0.25ex\right\vert\kern-0.25ex\right\vert}}

\newtheorem{theorem}{Theorem}
\theoremstyle{plain}

\newtheorem{definition}{Definition}
\newtheorem{lemma}{Lemma}
\newtheorem{proposition}{Proposition}
\newtheorem{remark}{Remark}

\numberwithin{equation}{section}
\numberwithin{theorem}{section}
\numberwithin{corollary}{section}
\numberwithin{definition}{section}
\numberwithin{lemma}{section}
\numberwithin{proposition}{section}
\numberwithin{remark}{section}

\title{A mixed formulation for the fractional Poisson problem}

\author{Juan Pablo Borthagaray}
\author{Nahuel de León}
\address{PEDECIBA -- Program for the Development of Basic Sciences and Facultad de Ingeniería, Universidad de la República, Av. Julio Herrera y Reissig 565, 11300 Montevideo, Departamento de Montevideo, Uruguay.}

\begin{document}

\begin{abstract} 
    The mixed formulation of the classical Poisson problem introduces the flux as an additional variable, leading to a system of coupled equations.
    Using fractional calculus identities, in this work we explore a mixed formulation of the fractional Poisson problem and establish its well-posedness. 
    Since a direct discretization of this problem appears to be out of reach, we adapt a stabilized approach by Hughes and Masud, which yields a coercive and well-posed formulation.
    The coercivity ensures the stability of any conforming finite element discretization. We further prove the convergence of this discretization, derive convergence rates, and discuss implementation aspects. Finally, we present numerical experiments that highlight both the importance of stabilization and the accuracy of our theoretical results.
\end{abstract}
\maketitle

\section{Introduction}
In recent years, study of nonlocal operators has been
an active area of research in different branches of mathematics. Nonlocal models have been increasingly used in different areas of science such as machine learning \cite{lu2022nonparametriclearningkernelsnonlocal,JMLR:v11:rosasco10a}, finance \cite{4fc30f5f-894d-3c6c-9b15-01f8f8d74820}, image processing \cite{buades2010image, gilboa2007nonlocal,lou2010image}, 
magnetohydrodynamics \cite{schekochihin2008mhd}, among others.
In particular, the fractional Laplacian has been considered in many applications, including, for example, diffusion-reaction problems \cite{yamamoto2012asymptotic}, quasi-geostrophic flows \cite{constantin1999behavior}, transport in porous media \cite{de2012general}, and ultrasound \cite{treeby2010modeling}.

Let $\Omega \subset \R^d$ be a bounded, Lipschitz domain. In this paper, we consider the fractional Poisson problem in $\Omega$, namely
\begin{equation} \label{eq:Poisson} \tag{P}
\left\lbrace \begin{aligned}
(-\Delta)^s u = f & \quad  \mbox{ in } \Omega, \\
u = 0 &  \quad \mbox{ in } \Omega^c := \R^d \setminus \Omega .
\end{aligned} \right.
\end{equation}
We propose to study the mixed formulation of this problem, cf. \eqref{eq:Darcy} below, to approximate both the solution $u$ and its so-called fractional gradient in primal form. 
Above, $s \in (0,1)$, $f \in L^2(\Omega)$, and $(-\Delta)^s$ denotes the fractional Laplacian
\begin{equation} \label{eq:def-FL}
(-\Delta)^s u (x) := \nu(d,s) \mbox{ p.v.} \int_{\R^d} \frac{u(x)-u(y)}{|x-y|^{d+2s}} dy, \quad x \in \R^d,
\end{equation}
where
\begin{equation} \label{eq:def-cds}
\nu(d,s) := \frac{2^{2s} s \Gamma(s+\frac{d}2)}{\pi^{d/2} \Gamma(1-s)}.
\end{equation}
Note that $(-\Delta)^s$ is an operator of order $2s$; for a summary of basic properties of this operator, we refer to \cite{daoud2022fractional, di2012hitchhiker's,lischke2019fractionallaplacian}.

We introduce some notation. Given two vectors $v,w\in\R^d$, their Euclidean inner product will be denoted by $v\cdot w$ and Euclidean norms will be denoted by $|v|:=\sqrt{v\cdot v}$. For $A, B > 0$, we write $A \simeq B$ to indicate that both inequalities $A \le C B$ and $B \le C A$ hold with constants independent of $A,B$.

There are several nonequivalent notions of nonlocal differential operators, which in turn give rise to different ways of writing \eqref{eq:Poisson} in divergence form. One customary option is to consider  {\em unweighted} operators \cite{Du_et_al}, namely to regard gradients as two-point operators and the divergence operator as minus the adjoint of the gradient. In our setting, this would reduce to consider, for sufficiently smooth functions $w \colon \R^d \to \R$, $\vPsi \colon \R^d \times \R^d \to \R^d$,
\[\begin{aligned}
\mathcal{G}^s w \colon  \R^d \times \R^d \to \R^d, & \quad
\mathcal{G}^s w(x,y) := 
\sqrt{\frac{\nu(d,s)}2} \, \frac{(w(x)-w(y))}{|x-y|^{d/2+s}}\, \frac{(x-y)}{|x-y|}, \\
\mathcal{D}^s \vPsi \colon \R^d \to \R, & \quad \mathcal{D}^s \vPsi(x) := \sqrt{\frac{\nu(d,s)}2} \, \int_{\R^d} \frac{(\vPsi(x,y) + \vPsi(y,x) )}{|x-y|^{d/2+s}} \cdot \frac{(x-y)}{|x-y|} \, dy .
\end{aligned}
\]
This yields the identity $(-\Delta)^s = -\mathcal{D}^s \circ \mathcal{G}^s$. Although these notions arise naturally in the nonlocal setting, interpreting the fractional gradient as a two-point operator complicates its physical interpretation and obscures the classical perspective of gradients as indicating the direction of maximal growth. Since the above definition essentially corresponds to a difference quotient, we argue that in applications it may be more meaningful to construct nonlocal gradients by  integration of difference quotients.

In this vein, the fractional Laplacian can be regarded as a composition of certain {\em weighted,} non-local, vector calculus operators. Namely, given $w \colon \R^d \to \R$ smooth enough, we define its {\em fractional gradient of order $s$}, $\grads w \colon \R^d \to \R^d$,
\begin{equation}\label{eq:def-grads}
\grads w (x) := 
\begin{cases}
    \mu(d,s) \int_{\R^d}  \frac{(w(y)-w(x))}{|x-y|^{d+s}} \, \frac{(y-x)}{|x-y|} \, dy&\text{if $s\in(0,1)$,}\\
    \nabla w(x) &\text{if $s=1$,}\\
    \mu(d,s)\int_{\R^d}\frac{(w(y)-w(x)-\nabla w(x)\cdot(x-y))}{|x-y|^{d+s}} \, \frac{(x-y)}{|x-y|} \, dy&\text{if $s\in(1,2)$.}
\end{cases}
\end{equation}
and, given $\vPsi \colon \R^d \to \R^d$ we define its  {\em fractional divergence of order $s$}, $\divs \vPsi \colon \R^d \to \R$,
\begin{equation}\label{eq:def-divs}
\divs \vPsi (x) := \begin{cases}
    \mu(d,s) \int_{\R^d}  \frac{(\vPsi(y)-\vPsi(x))}{|x-y|^{d+s}} \cdot \frac{(y-x)}{|x-y|} \, dy&\text{if $s\in(0,1)$,}\\
    \div\, \vPsi(x) &\text{if $s=1$,}\\
    \mu(d,s)\int_{\R^d}\frac{(\vPsi(y)-\vPsi(x)-\nabla \vPsi(x)(x-y))}{|x-y|^{d+s}} \cdot \frac{(x-y)}{|x-y|} \, dy&\text{if $s\in(1,2)$.}
\end{cases}
\end{equation}
In the two definitions above, we have taken
\begin{equation}
\mu(d,s):=\frac{2^s\Gamma(\frac{d+s+1}{2})}{\pi^{d/2}\Gamma(\frac{1-s}{2})}.
\end{equation}
The operators $\grads$ were introduced by Horv\'ath \cite{Horvath} as generalizations of the Riesz transform, and interpreted as fractional gradients by Shieh and Spector \cite{Shieh2015}.
In \cite{vsilhavy2020fractional}, 
Šilhavý discuses the feasibility of \eqref{eq:def-grads} and later \eqref{eq:def-divs} as alternatives to the coordinate-base definitions for fractional derivatives. In that reference, it was shown that in the space of rapidly decaying Schwartz functions, the $\grads$ operator is uniquely determined by the properties of being rotationally invariant\footnote{Namely, that for every $Q\in O(d)$ it holds that $\grads[w(Q^T\,\cdot)](x)=Q\,\grads w(Q^Tx)$ (resp. $\divs[Q\vPsi(Q^T\,\cdot)](x)=\divs\vPsi(Q^Tx)$).}, translationally invariant, and $s$-homogeneous\footnote{Namely, that for all $\lambda>0$ it holds that $\grads [w(\lambda\,\cdot)](x)=\lambda^s\,\grads w(\lambda x)$ (resp. $\divs [\vPsi(\lambda\,\cdot)](x)=\lambda^s\,\divs \vPsi(\lambda x)$).}, i.e., that every operator possessing these properties is a multiple of $\grads$. An analogous result is valid for $\divs$; we refer to \cite{vsilhavy2020fractional} for details.
Furthermore, while we have provided the definition of these operators for orders $s \in (0,2)$, they can be defined for arbitrary $s \in \mathbb{C}$ satisfying $\mbox{Re}(s) > -d$. In particular, in this work we make use of \eqref{eq:def-grads} for $s \in (0,2)$ but we will only need \eqref{eq:def-divs} for $s \in (0,1)$.

We make a brief discussion on some properties of these operators, and refer to \cite{Comi2019, Comi2023,delia2021unifiedtheoryfractionalnonlocal,vsilhavy2020fractional} for more details. In first place, one can resort to the Fourier transform to express, for functions in the Schwartz class, that \cite[Theorem 1.4]{Shieh2015}
\begin{equation}\label{eq:Fourier-oper} \begin{split}
 \widehat{(-\Delta)^s w}(\xi) &= (2\pi)^{2s}|\xi|^{2s} \widehat{w}(\xi), \\
 \widehat{\grads w}(\xi) & = (2\pi)^s i |\xi|^{s-1} \xi \widehat{w}(\xi), \\
 \widehat{\divs \Psi}(\xi) & = (2\pi)^s i |\xi|^{s-1} \xi \cdot \widehat{\Psi}(\xi),
\end{split}\end{equation}
from which the formula 
\begin{equation}\label{eq:frac-comp}(-\Delta)^s w = -\divs \grads w,\end{equation}
follows immediately\footnote{Actually, upon a suitable definition of the operators, the identity $(-\Delta)^s w = -\div^\alpha \grad^\beta w$ holds whenever $\mbox{Re }\alpha > -d$, $\mbox{Re }\beta \ge 0$, and $s = \frac{\alpha+\beta}2$, see \cite[Theorem 5.3]{vsilhavy2020fractional}.}
(see also \cite[Proposition 5.3]{delia2021connectionsnonlocaloperatorsvector} and \cite[Theorem 5.3]{vsilhavy2020fractional} for different approaches). 
Identity \eqref{eq:frac-comp} is analogous to the classical case and serves as the main motivation for this work.
Identities \eqref{eq:Fourier-oper} allow one to rewrite the fractional-order operators as a convolution of their classical counterparts with the well-known Riesz kernels \cite[Chapter 5]{Stein}. Indeed, for $0<\alpha<d$ and $x\not=0$, the Riesz kernel of order $\alpha$ is defined as 
\[I_\alpha(x)=c_{\alpha,d}|x|^{-d+\alpha},\quad c_{\alpha,d}=\frac{\Gamma((d-\alpha)/2)}{\pi^{d/2}2^\alpha\Gamma(\alpha/2)}.\]
We have that $\widehat{I_\alpha}(\xi)=|2\pi\xi|^{-\alpha}$, so we deduce
\begin{equation}\label{eq:grads-divs-riesz-def}  \grads w =I_{1-s}*\nabla w, \qquad
 \divs \Psi =I_{1-s}*\div\, \Psi,
\end{equation}
for Lipschitz functions with compact support  \cite[Section 2.4]{Comi2019}.

Additionally, we have an integration by parts formula (e.g. \cite[Section 6]{vsilhavy2020fractional}): given $w \in C^\infty_c(\R^d)$, $\vPsi \in C^\infty_c(\R^d, \R^d)$,
\begin{equation} \label{eq:ibp}
\int_{\R^d} \grads w \cdot \vPsi = - \int_{\R^d} w \, \divs \vPsi.
\end{equation}
Using a density argument, we can extend this formula to a broader class of functions, namely the spaces $\tHs$ and $H(\divs; \Omega)$ defined below.

We consider the well-known fractional Sobolev space
\[H^s(\R^d):=\left\{w\in L^2(\R^d):|w|_{H^s(\R^d)}<\infty\right\},\]
 where
\begin{equation} \label{eq:def-norma-Hs}
|w|_{H^s(\R^d)}:=\left(\frac{\nu(d,s)}{2}\int_{\R^d}\int_{\R^d}\frac{|w(x)-w(y)|^2}{|x-y|^{d+2s}}dx\,dy\right)^{\frac12}
=\| \grads w \|_{L^2(\R^d)}.
\end{equation}
We refer to \cite{di2012hitchhiker's} for an introduction to the topic and the main properties of this space, and we also define
\begin{equation}\label{eq:def-tHs}
\tHs := \{ w \in H^s(\R^d) \colon \mbox{supp } w \subset \overline\Omega \}.
\end{equation}
We have the following Poincar\'e inequality (see, for example, \cite[Theorem 3.9]{edmunds2022fractional}),
\begin{equation}\label{eq:Poincare}
\| w \|_{L^2(\Omega)} \le C_P |w|_{H^s(\R^d)} 
, \quad \forall w \in \tHs,
\end{equation}
that allows us to consider the norm 
\[
\| w \|_{\tHs} := |w|_{H^s(\R^d)}.
\]

We will also make use of the space
\begin{equation}\label{eq:def-Hdivs}
H(\divs; \Omega) := \Bigg\{ \vPsi \in \Big(L^2(\R^d)\Big)^d \colon (\divs \vPsi) \big|_\Omega \in L^2(\Omega) \Bigg\},
\end{equation}
with the norm
\[
\| \vPsi \|_{H(\divs; \Omega)} = \left( \| \vPsi \|_{L^2(\R^d)}^2 + \| (\divs \vPsi) \big|_\Omega \|_{L^2(\Omega)}^2 \right)^{1/2}.
\]
Furthermore, we will denote by $\widetilde{L}^2(\Omega)$ the space of functions in $L^2(\Omega)$ that are extended by zero to $\Omega^c$.

In this work, we analyze the mixed formulation of the fractional Poisson problem, to which we will refer as the {\em fractional Darcy problem}.
With the notation defined above, that problem reads: find $(p, \vPhi) \in \widetilde{L}^2(\Omega)\times H(\divs; \Omega)$ such that
\begin{equation} \label{eq:Darcy} \tag{D}
\left\lbrace \begin{aligned}
 \vPhi + \grads p = 0 & \quad  \mbox{ in } \R^d, \\
\divs \vPhi = f & \quad  \mbox{ in } \Omega, \\
p = 0 &  \quad \mbox{ in } \Omega^c .
\end{aligned} \right.
\end{equation}

Following the standard terminology for the Darcy problem, occasionally we will refer to $p$ as {\em pressure} and to $\vPhi$ as {\em flux}. Clearly, if $(p,\vPhi)$ solves the problem above, then $u=p$ solves \eqref{eq:Poisson}. The difference between \eqref{eq:Poisson} and \eqref{eq:Darcy} is that, obviously, we have introduced the flux variable $ \vPhi  := - \grads p$; due to the nonlocal nature of the problem, this definition needs to be imposed in the whole space $\R^d$ and not just in the domain $\Omega$. 

From a computational perspective, in comparison to the classical (local) divergence, the operator $\divs$ brings two apparent difficulties. On the one hand, it does not map piecewise polynomial functions into piecewise polynomials of one degree less. On the other hand, the nonlocal nature of the operator implies that $\divs \vPsi$ can have an unbounded support even when $\vPsi$ does.
Thus, the construction of $H(\divs)$-conforming finite elements a-la Raviart-Thomas or Brezzi-Douglas-Marini seems unfeasible. For this reason, in this work we adapt the approach of Masud and Hughes \cite{Masud2002} and employ a stabilized formulation that can be discretized with continuous Lagrange elements.

The structure of this work is as follows. In Section \ref{sec:problem-formulation}, we establish the well-posedness of problem \eqref{eq:Darcy}. Section \ref{sec:stabilization} introduces a stabilized formulation: we show it to be continuous and coercive in a suitable norm, thereby obtaining a well-posed problem. Section \ref{sec:regularity} deals with the Sobolev regularity of solutions to our problem. We show estimates for both variables being approximated, including explicit decay bounds for the flux and weighted estimates for fractional gradients of the pressure. The finite element framework we employ is described in Section \ref{sec:FE}, where some technical results on interpolation are derived and a priori error bounds are shown. This work concludes with several numerical experiments, displayed in Section \ref{sec:experiments}, that illustrate our theoretical findings.

\section{Mixed formulation and well-posedness} \label{sec:problem-formulation}

This section analyzes the well-posedness of the fractional Darcy problem \eqref{eq:Darcy}. For simplicity, we assume the right-hand side $f \in L^2(\Omega)$.
Using the integration by parts formula \eqref{eq:ibp}, its weak formulation reads: find $(p, \vPhi) \in \widetilde{L}^2(\Omega)\times H(\divs; \Omega)$ such that, for all $(q,\vPsi) \in \widetilde{L}^2(\Omega)\times H(\divs; \Omega)$,
\begin{equation} \label{eq:Darcy-weak}
\int_{\R^d} \vPhi \cdot \vPsi - \int_{\R^d} p \, \divs \vPsi +
\int_{\R^d} q \, \divs \vPhi = \int_{\R^d} f q.
\end{equation}
We point out that all but the first of the integrals above can be actually computed in $\Omega$.

Let us introduce some additional notation. First, we define the forms 
\begin{equation}\label{eq:def-ab} \begin{aligned}
a \colon H(\divs; \Omega) \times H(\divs; \Omega) \to \R, & \quad a (\vPhi, \vPsi) = \int_{\R^d} \vPhi \cdot \vPsi, \\
b \colon \widetilde{L}^2(\Omega)\times H(\divs; \Omega) \to \R, & \quad b(q,\vPsi) = \int_{\Omega} q \, \divs \vPsi , \\
F \colon \widetilde{L}^2(\Omega) \to \R, & \quad F(q) =  \int_{\Omega} f q .
\end{aligned} \end{equation}
We also define
$
B \colon H(\divs; \Omega) \to L^2(\Omega),$
\[
(B \vPsi , q)_{L^2(\Omega)} := b(q, \vPsi) , \quad \forall \vPsi \in H(\divs; \Omega), \ q \in \widetilde{L}^2(\Omega).
\]
Notice that $a$ is symmetric and that the problem above has a clear saddle-point structure. Therefore, by standard arguments in the analysis of mixed formulations, to prove the well-posedness of \eqref{eq:Darcy-weak} it suffices to show that
\begin{itemize}
\item the form $a$ is coercive in $\ker B$;
\item the form $b$ satisfies an inf-sup condition.
\end{itemize}

The fact that $a$ is coercive in $\ker B$ follows straightforwardly upon observing that $\ker B = \{ \vPsi \in H(\divs; \Omega) \colon \divs \vPsi = 0 \mbox{ in } \Omega \}$. This yields that, for every $\vPsi \in \ker B$,
\[
a(\vPsi, \vPsi) = \| \vPsi \|_{L^2(\R^d)}^2 =  \| \vPsi \|_{H(\divs; \Omega)}^2.
\]

\begin{lemma}[surjectivity of $\divs$]
Let $\Omega \subset \R^d$ be a bounded, Lipschitz domain. The operator $\divs |_\Omega$ such that $\divs |_\Omega \vPsi := (\divs \vPsi)|_\Omega$ maps $H(\divs; \Omega)$ onto 
$\widetilde{L}^2(\Omega)$. Consequently, $b$ satisfies an inf-sup condition: there exists $\beta > 0$ such that
\begin{equation}\label{eq:inf-sup-b}
\inf_{p \in \widetilde{L}^2(\Omega)} \sup_{\vPhi \in H(\divs; \Omega)} \frac{b(p,\vPhi)}{\| p \|_{L^2(\Omega)} \|\vPhi\|_{H(\divs; \Omega)}} \ge \beta .
\end{equation}
\end{lemma}
\begin{proof}
The fact that $\divs |_\Omega \vPsi \in L^2(\Omega)$ for all $\vPsi \in H(\divs; \Omega)$ is evident from the definition. \
Next, given $p \in \widetilde{L}^2(\Omega)$, we seek $w \in \tHs$ such that
\[
(\grads w, \grads v)_{L^2(\R^d)} = (p,v)_{L^2(\R^d)}.
\]
By the identity \eqref{eq:frac-comp} and the integration by parts formula \eqref{eq:ibp}, this is equivalent to stating that $w \in \tHs$ solves the fractional Poisson problem \eqref{eq:Poisson} with right-hand side $p$. We let $\vPsi := - \grads w$, that obviously satisfies $\divs \vPsi = (-\Delta)^s w = p$ in $\Omega$. This shows that $\divs |_\Omega$ is surjective.
Additionally, we have
\[
\| \vPsi \|_{L^2(\R^d)}^2 = \| \grads w \|_{L^2(\R^d)}^2 = (p,w)_{L^2(\R^d)} \le \| p \|_{L^2(\Omega)} \| w \|_{L^2(\Omega)},
\]
 and the Poincar\'e inequality \eqref{eq:Poincare} gives $$ \| w \|_{L^2(\Omega)} \le C_P \|\grads w \|_{L^2(\R^d)} = C_P \| \vPsi\|_{L^2(\R^d)},$$ so that $\| \vPsi \|_{L^2(\R^d)} \le C_P\| p \|_{L^2(\Omega)}$.
 We therefore deduce that
 \[ \begin{aligned}
 \| \vPsi \|_{H(\divs; \Omega)}^2 & =  \| \vPsi \|_{L^2(\R^d)}^2 +  \| \divs \vPsi |_\Omega \|_{L^2(\Omega)}^2  \le (1 + C_P^2) \| p \|_{L^2(\Omega)} ^2.
 \end{aligned} \]
 
 Finally, we conclude
 \[
 \sup_{\vPhi \in H(\divs; \Omega)} \frac{b(p,\vPhi)}{\|\vPhi\|_{H(\divs; \Omega)}}  \ge  \frac{b(p,\vPsi)}{\|\vPsi\|_{H(\divs; \Omega)}} \ge  \frac{\| p \|_{L^2(\Omega)}}{\sqrt{1+C_P^2}}.
 \]
This shows that $b$ satisfies the inf-sup condition \eqref{eq:inf-sup-b} with constant $\beta := \sqrt{1+C_P^2}^{-1}$.
\end{proof}
As a corollary, using standard tools in the analysis of mixed formulations (see, for example, \cite[Theorem 4.2.3]{Boffi2013}), we deduce the well-posedness of the fractional Darcy problem.

\begin{proposition}[well-posedness]
Let $\Omega \subset \R^d$ be a bounded, Lipschitz domain. Then, problem \eqref{eq:Darcy-weak} has a unique solution $(p, \vPhi) \in \widetilde{L}^2(\Omega)\times H(\divs; \Omega)$, and there hold
\begin{equation} \label{eq:Darcy-continuity} \begin{aligned}
\| p \|_{L^2(\Omega)} & \le 2 \sqrt{1 + C_P^2}\,\| f \|_{L^2(\Omega)} , \\
\| \vPhi \|_{H(\divs; \Omega)} & \le 2 (1 + C_P^2)\, \| f \|_{L^2(\Omega)},
\end{aligned} \end{equation}
with $C_P$ being the constant in the Poincar\'e inequality \eqref{eq:Poincare}.
\end{proposition}

\section{Stabilization} \label{sec:stabilization}
Although the formulation developed in the previous section is well-posed, its saddle-point structure makes finite element approximation challenging. In this section, we build on the ideas of \cite{Masud2002} to develop a stabilized variational formulation of the fractional Darcy problem, which is amenable to be treated with continuous Lagrange elements. We also present numerical experiments that illustrate the importance of the stabilization.

\subsection{Stabilized problem}
To shorten the notation, we write 
\begin{equation}\label{eq:def-calL} \begin{aligned}
&\calL \colon \left(\widetilde{L}^2(\Omega)\times H(\divs; \Omega) \right)\times \left(\widetilde{L}^2(\Omega)\times H(\divs; \Omega) \right)\to \R, \\
& \calL ((p, \vPhi), (q, \vPsi)) := a(\vPhi,\vPsi) - b(p, \vPsi) + b(q, \vPhi) ,
\end{aligned}\end{equation}
so that we can rewrite \eqref{eq:Darcy-weak} as 
\begin{equation} \label{eq:Darcy-weak2}
\calL ((p, \vPhi), (q, \vPsi)) = F(q) .
\end{equation}
We introduce the stabilized form in $\bV :=  \tHs \times H(\divs; \Omega)$, $\calL_{\rm stab} \colon \bV \times \bV \to \R$, 
\begin{equation}\label{eq:def-calLstab2} \begin{aligned}
& \calL_{\rm stab} ((p, \vPhi), (q, \vPsi)) := \calL ((p, \vPhi), (q, \vPsi)) + \frac12 \int_{\R^d} \left(\vPhi + \grads p\right) \cdot \left(-\vPsi + \grads q\right) .
\end{aligned}\end{equation}

By \eqref{eq:def-calL}, we note that this can be rewritten as
\begin{equation}\label{eq:def-calLstab}
    \begin{aligned}
\calL_{\rm stab}((p, \vPhi), (q, \vPsi)) =& \frac12\int_{\R^d} \vPhi \cdot \vPsi+\frac12\int_{\R^d} \grads p \cdot \vPsi-\frac12\int_{\R^d} \grads q \cdot 
 \vPhi\\&+\frac12 \int_{\R^d} \grads p \cdot \grads q .
\end{aligned}
\end{equation}

With this, we consider the stabilized problem: find $(p, \vPhi) \in \bV$ such that
\begin{equation} \label{eq:Darcy-stabilized}
 \calL_{\rm stab} ((p, \vPhi), (q, \vPsi)) = F(q) \quad  \forall  (q, \vPsi) \in \bV.
\end{equation}

We make three important remarks regarding the definition of $\calL_{\rm stab}$. First, we have shrunk the domain of $\calL$ by replacing $\widetilde{L}^2(\Omega)$ by $\tHs$ so that the stabilization term is well-defined. Second, the stabilization term above needs to be computed in the whole $\R^d$. Third, it is consistent: in the solution of our problem, the term we are adding is zero.
Let us consider the following norm in $\bV$,
\begin{equation}\label{eq:def-vertiii}
\vertiii{(q, \vPsi)} := \left[ \frac12 \left( \| \grads q \|_{L^2(\R^d)}^2 + \| \vPsi \|_{L^2(\R^d)}^2 \right)\right]^{1/2}.
\end{equation}
The fact that $\vertiii{\cdot}$ is a seminorm is evident. We further 
note that if $\vertiii{(q, \vPsi)}=0$ then obviously $\vPsi=0$ in $\R^d$ and the Poincaré inequality \eqref{eq:Poincare} implies that $q=0$ in $\Omega$ (and thus in $\R^d$), and it follows that $\vertiii{\cdot}$ is actually a norm.

\begin{lemma}[equivalence]
A pair $(p, \vPhi) \in \bV$ solves the stabilized problem \eqref{eq:Darcy-stabilized} if and only if it solves \eqref{eq:Darcy-weak2}.
\end{lemma}

\begin{proof}
If $(p, \vPhi) \in \bV$ solves problem \eqref{eq:Darcy-weak2} then using $q=0$ in that formulation and integrating by parts (cf. \eqref{eq:ibp}), we deduce $\vPhi + \grads p = 0$ in $\R^d$. Therefore, the stabilization term vanishes,
\[
\frac12 \int_{\R^d} \left(\vPhi + \grads p\right) \cdot \left(-\vPsi + \grads q\right) = 0 \quad \forall (q,\vPsi) \in \bV,
\]
and it follows that $(p, \vPhi)$ solves the stabilized problem \eqref{eq:Darcy-stabilized}.

Conversely, if $(p, \vPhi) \in \bV$ solves the stabilized problem \eqref{eq:Darcy-stabilized}, then for every $\vPsi\in H(\divs;\Omega)$ we have
\[
\begin{aligned}
    0=\calL_{\rm stab}((p,\vPhi),(0,\vPsi))&=\frac12\int_{\R^d}\vPhi \cdot \vPsi-\frac12\int_{\R^d}p\,\divs\vPsi\\
    &=
    \frac12\int_{\R^d}\vPhi\cdot\vPsi+\frac12\int_{\R^d}\grads p\cdot\vPsi,
    \end{aligned}
\]
which again implies $\grads p + \vPhi = 0$ in $\R^d$. Therefore, we deduce
\[
\calL((p,\vPhi)(q,\vPsi))=\calL_{\rm stab}((p,\vPhi)(q,\vPsi))=F(q).
\]
\end{proof}

We have obtained an equivalent formulation to the fractional Darcy problem, but the stabilized form has the key advantage of being coercive.

\begin{lemma}[coercivity]
We have 
\[
\calL_{\rm stab}((p, \vPhi), (p, \vPhi)) = \vertiii{(p, \vPhi)}^2 \quad \forall (p, \vPhi) \in  \bV.
\]
\end{lemma}
\begin{proof}
The result follows by a direct computation using \eqref{eq:def-calLstab}. Indeed, if $(p, \vPhi) \in  \bV$, then
\begin{align*}
\calL_{\rm stab}((p, \vPhi), (p, \vPhi)) &
=\frac{1}{2}\|\vPhi\|_{L^2(\R^d)}^2+\frac{1}{2}\|\grads p\|_{L^2(\R^d)}^2.
\end{align*}
\end{proof}

In addition to being coercive, it is straightforward to verify that the form $\calLs$ is continuous.

\begin{lemma}[continuity]
We have 
\[
\calL_{\rm stab}((p, \vPhi), (q, \vPsi)) \le  \vertiii{(p, \vPhi)}  \, \vertiii{(q, \vPsi)} \quad \forall (p, \vPhi), (q, \vPsi) \in \bV.
\]
\end{lemma}
\begin{proof}
Let $(p, \vPhi), (q, \vPsi) \in \bV$. Using \eqref{eq:def-calLstab}, we have 
\begin{equation}
\begin{aligned}
|\calL_{\rm stab}((p, \vPhi), (q, \vPsi))| =& \bigg|\frac12\int_{\R^d} \vPhi \cdot \vPsi+\frac12\int_{\R^d} \grads p \cdot \vPsi-\frac12\int_{\R^d} \grads q \cdot 
 \vPhi\\&+\frac12 \int_{\R^d} \grads p \cdot \grads q \bigg|.
\end{aligned}
\end{equation}
Therefore, by Young's inequality, we deduce
\begin{equation}
\begin{aligned}
|\calL_{\rm stab}((p, \vPhi), (q, \vPsi))| &\le \frac12\left(\|\vPhi\|_{L^2(\R^d)}+\|\grads p\|_{L^2(\R^d)}\right)\left(\|\vPsi\|_{L^2(\R^d)}+\|\grads q\|_{L^2(\R^d)}\right)\\
&\le \,\vertiii{(p,\vPhi)} \, \vertiii{(q,\vPsi)}.
\end{aligned}
\end{equation}
\end{proof}

Finally, the combination of the two lemmas above with the Lax-Milgram theorem gives rise to the well-posedness of our problem.

\begin{proposition}[well-posedness of stabilized formulation]\label{prop:stability}
Given $f \in L^2(\Omega)$, problem \eqref{eq:Darcy-stabilized} has a unique solution $ (p, \vPhi) \in \bV$. Moreover, we have the stability estimate
\[
\vertiii{(p, \vPhi)} \le \sqrt{2} \, C_P \| f \|_{L^2(\Omega)},
\]
with $C_P$ being the constant in the Poincar\'e inequality \eqref{eq:Poincare}.
\end{proposition}
\begin{proof}
The Lax-Milgram theorem implies the existence and uniqueness of a solution $(p, \vPhi) \in \bV$. 
Additionally, by the coercivity of the stabilized form and identity \eqref{eq:Darcy-stabilized}, we deduce 
\[
\vertiii{(p, \vPhi)}^2 = \calL_{\rm stab}((p, \vPhi),(p, \vPhi)) = F(p) \le \| p \|_{L^2(\Omega)} \| f \|_{L^2(\Omega)}.
\]
The desired estimate follows now by the Poincaré inequality \eqref{eq:Poincare} and noticing that $\| p\|_{L^2(\R^d)} \le \sqrt{2}\vertiii{(p, \vPhi)}$.
\end{proof}

\subsection{The need of stabilization}
Our aim in introducing the stabilized formulation is to be able to use standard $\mathcal{P}_1$-$\mathcal{P}_1^d$ discretization for problem \eqref{eq:Darcy}. As in the local case, the stability of the standard mixed formulation of the fractional Laplacian is not guaranteed when continuous Lagrange elements are employed for both the pressure and the flux.

We illustrate this point with a computational experiment. Using the finite element approximations described in Section \ref{sec:FE}, and implementing the matrices as outlined in Section \ref{sec:experiments}, we consider problem \eqref{eq:Darcy} with constant right-hand side $f=1$ on the square domain $\Omega = (-1,1)^2$.

\begin{figure}[!htbp]
    \centering
    \begin{minipage}{0.4\textwidth}
    \includegraphics[width=\textwidth]{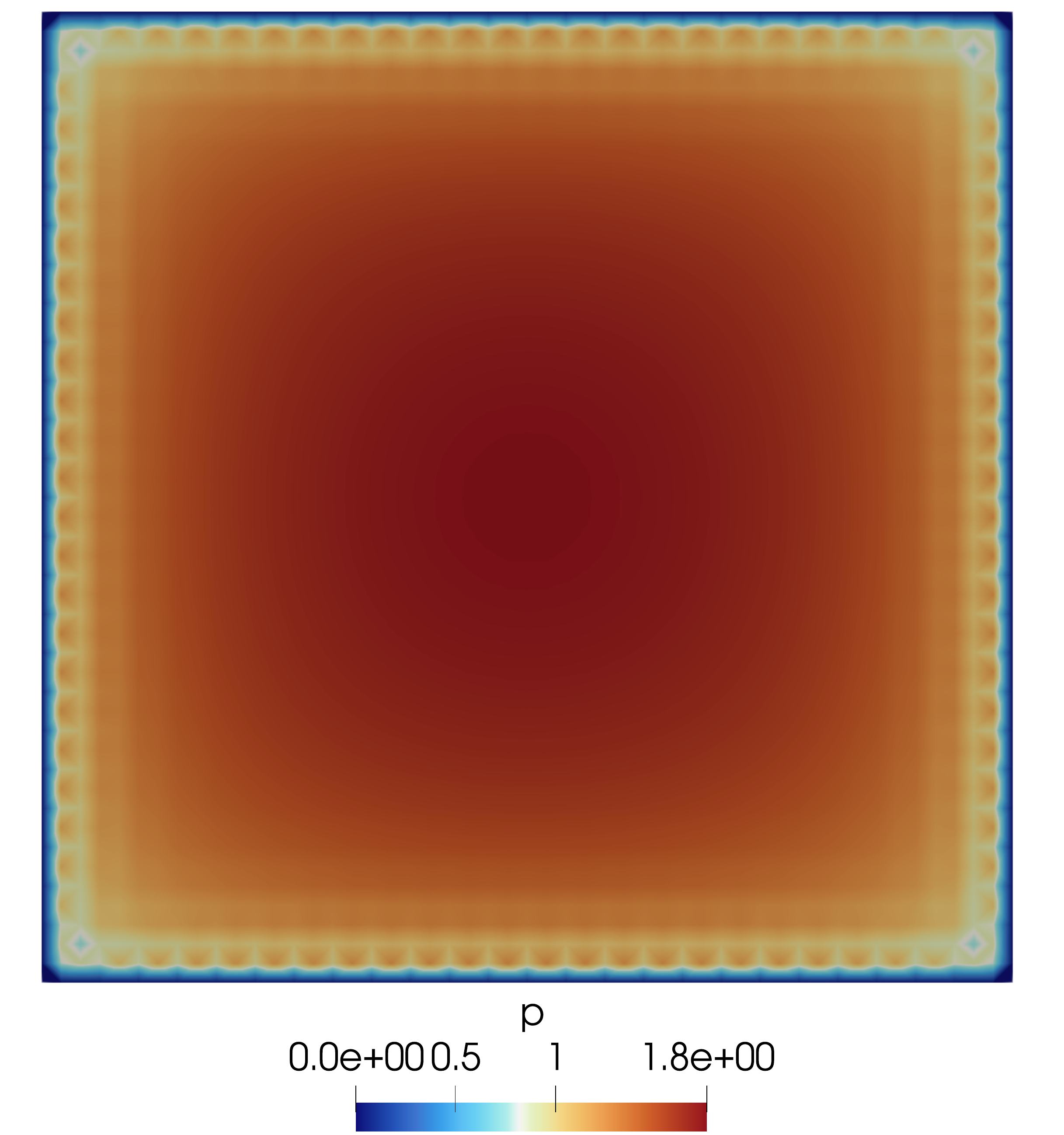}
    \end{minipage}
    \hfill
    \begin{minipage}{0.4\textwidth}
    \includegraphics[width=\textwidth]{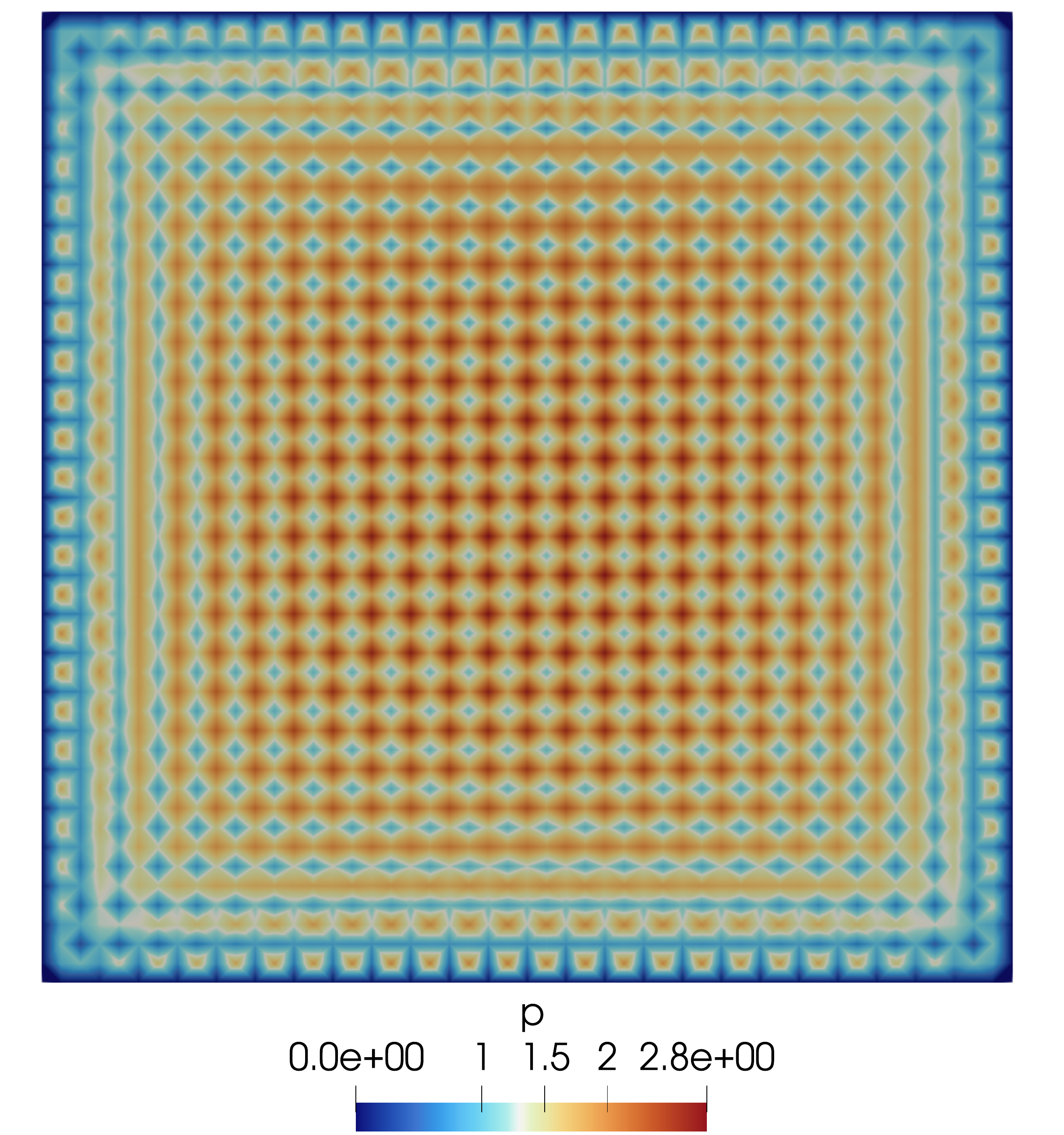}
    \end{minipage}
    \begin{minipage}{0.4\textwidth}
    \includegraphics[width=\textwidth]{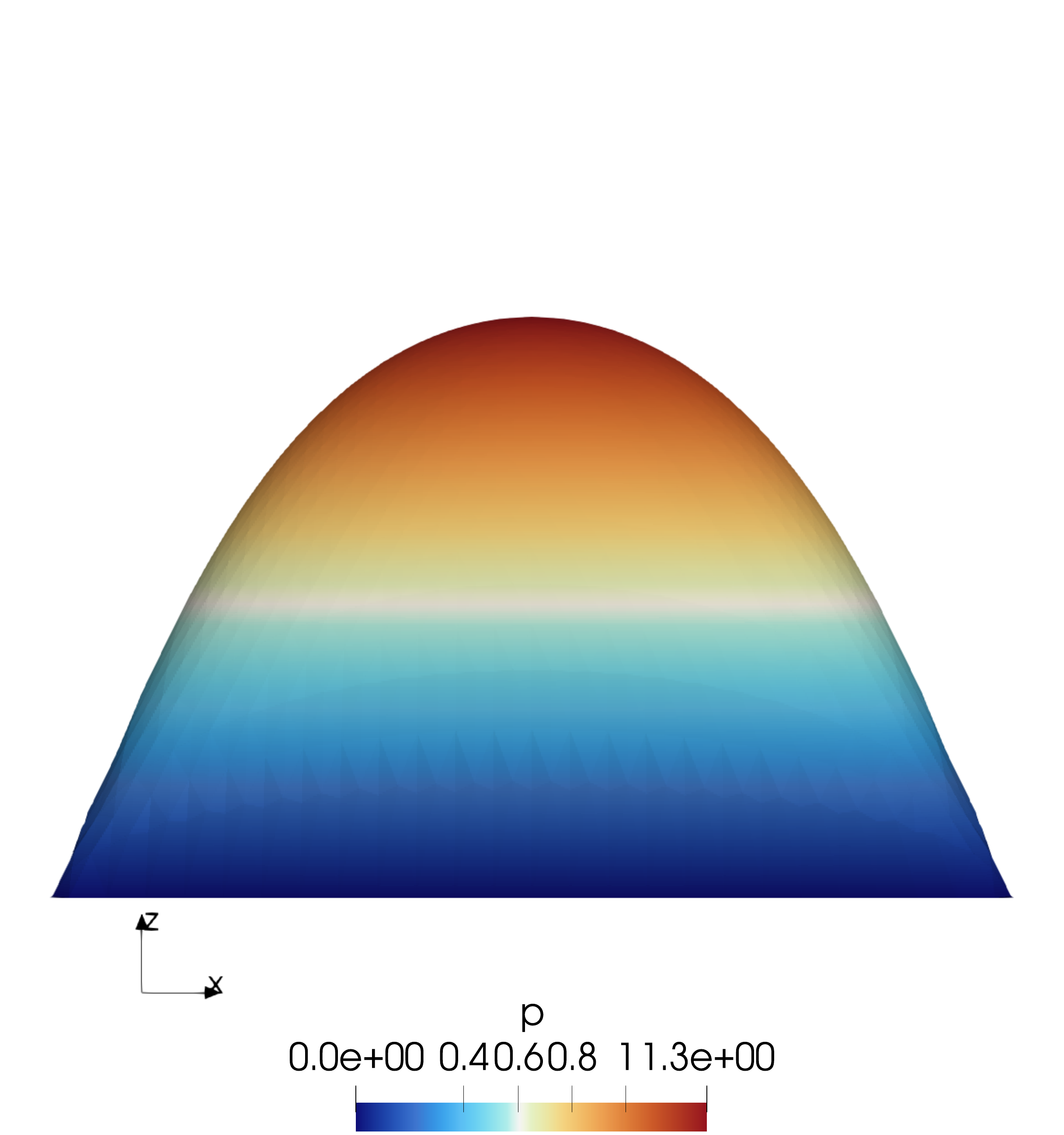}
    \end{minipage}
    \hfill
    \begin{minipage}{0.4\textwidth}
    \includegraphics[width=\textwidth]{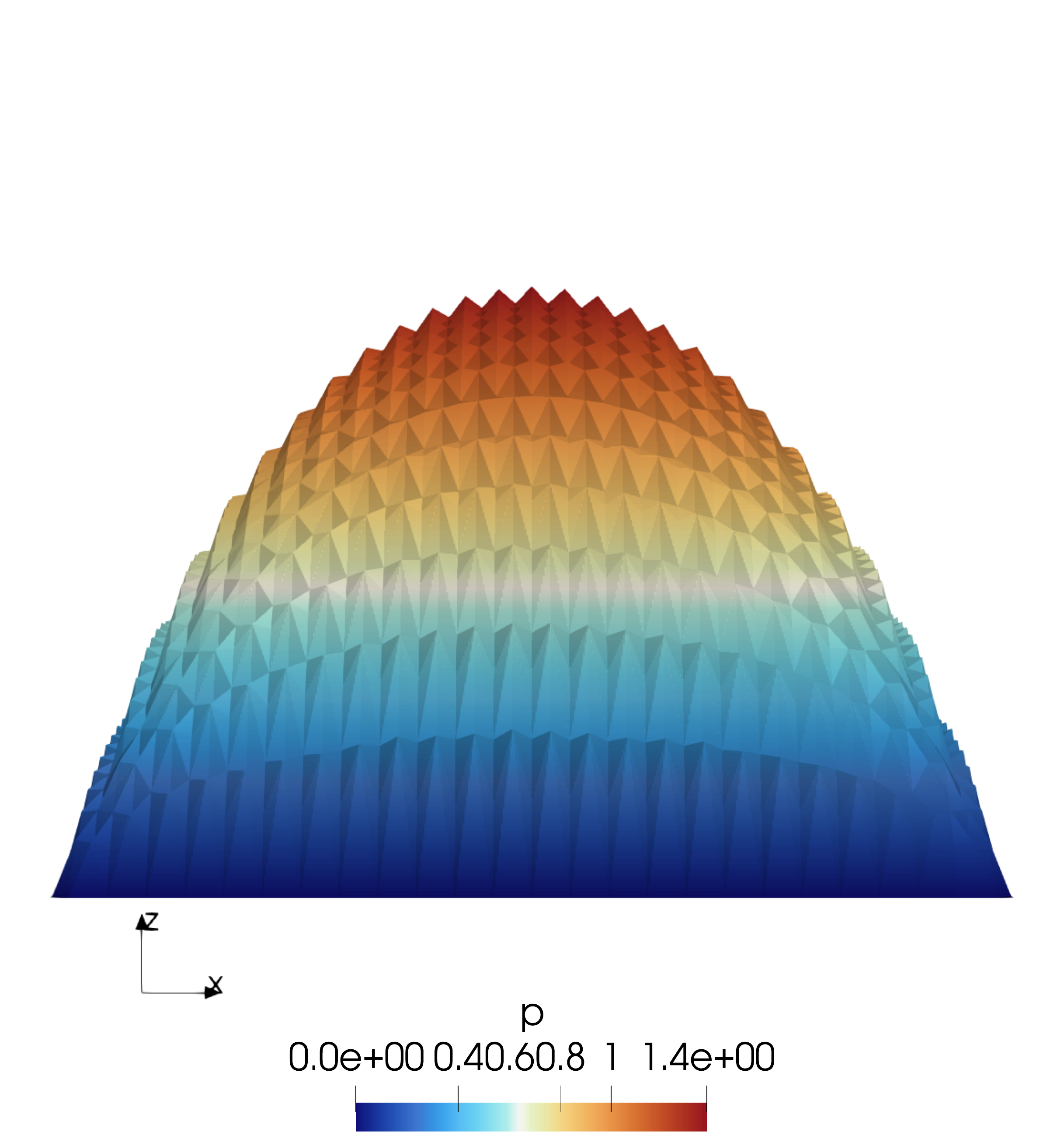}
    \end{minipage}
    \caption{Finite element approximations of the pressure $p$ in \eqref{eq:Darcy} for $s=0.25$ (top row) and $s=0.75$ (bottom) with the stabilization factor  (left column) and without it (right).    
    Here, $f\equiv 1$, $\Omega=(-1,1)^2$, the mesh size is $h=0.04$ and the computational domain $B_H=B(0,2)$. See Section \ref{sec:FE} for details on the finite element framework.}  \label{fig:inestable_vs_estable}
\end{figure}

Figure \ref{fig:inestable_vs_estable} shows some pressures computed  with continuous, piecewise linear finite elements on structured meshes with $s=0.25$ and $s=0.75$ for the finite element counterparts of \eqref{eq:Darcy-weak} and \eqref{eq:Darcy-stabilized}. The results clearly show that the non-stabilized formulation produces spurious oscillations.

To further illustrate this observation, we compute the pressure $H^s$-errors on quasi-uniform meshes for problem \eqref{eq:exp1} --for which an explicit solution is available-- with $s=0.5$ for both the stabilized and the non-stabilized formulations. The results show that the errors for the standard mixed formulation are significantly larger than those for its stabilized counterpart.

\begin{table}[!htbp]
    \centering
    \begin{tabular}{c|c|c|c|c} 
         & $h=0.1$ & $h=0.05$ & $h=0.025$ & $h=0.02$ \\
        \hline
Non-stabilized & 0.7908 & 0.5833 & 0.4338 & 0.4124 \\ \hline
Stabilized & 0.1056 & 0.0705 & 0.0488 & 0.0443 \\ 
    \end{tabular}
    \caption{$H^s$-error for the pressure $p$ on problem \eqref{eq:exp1} with $s=0.5$, using the standard (without stabilization) and stabilized mixed formulations, respectively.}
\end{table}

\section{Sobolev regularity} \label{sec:regularity}

Coercivity ensures that any conforming discretization satisfies a best-approximation property. To derive convergence rates, however, interpolation estimates and solution regularity are required. We address the latter in this section.

Sobolev regularity up to $\partial\Omega$ of $u$, the solution of \eqref{eq:Poisson}, was established in \cite{Borthagaray2023}; for $f$ in the Besov space $B^{-s+1/2}_{2,1}(\Omega)$, the solution $u$ lies in the space $\cap_{\eps > 0} \widetilde{H}^{s+1/2-\eps}(\Omega)$. This is the maximal expected regularity not only for arbitrary Lipschitz domains, but also for smooth domains as well. Indeed, a remarkable example arises when $\Omega=B(0,1)$ and $f\equiv1$. The solution in this case corresponds to the first exit time of a $2s$-stable L\'evy process in $\Omega$. The solution is given by 
\[u(x)=k_{s,d}(1-|x|^2)^s\chi_{B(0,1)}(x),\]
with $k_{s,d}>0$ an explicit constant (cf. Theorem \ref{thm:explicit_solutions}). Despite both the domain and the right hand side being smooth, the solution satisfies
\[
u\in\cap_{\eps > 0} \widetilde{H}^{s+1/2-\eps}(\Omega),\quad u\notin\widetilde{H}^{s+1/2}(\Omega).
\]

We point out that the hypothesis $f\in B^{-s+1/2}_{2,1}(\Omega)$ is weaker than $f \in L^2(\Omega)$ when $s > 1/2$; if $s\le 1/2$, one can perform a simple interpolation argument to show that $u \in \cap_{\eps > 0} \widetilde{H}^{2s-\eps}(\Omega)$ provided $f \in L^2(\Omega)$. 

Moreover, having at hand the regularity of $u$, one can deduce a regularity estimate for the flux by means of standard mapping properties of $\grads$. We include a proof of such properties for completeness.
\begin{lemma}[continuity of $\grads$]
    \label{lemma:grads-mapping}
    Let $r,s\geq0$. Then, for all $u\in H^{r+s}(\R^d)$ it holds
    \[
        \|\grads u\|_{H^r(\R^d)}\leq 2\pi\|u\|_{H^{r+s}(\R^d)}.
    \]
    In particular, if $u\in \widetilde{H}^{r+s}(\Omega)$ then $\grads u\in ({H}^{r}(\R^d))^d$.
\end{lemma}
\begin{proof}
 Given $u\in H^{r+s}(\R^d)$, we apply the identities \eqref{eq:Fourier-oper} and \eqref{eq:def-norma-Hs} to obtain
    \[
    \begin{aligned}
        |\grads u|_{H^r(\R^d)}&=\|\grad^r(\grads u)\|_{L^2(\R^d)}\\
        &=
        \|(2\pi)^r|\cdot|^{r}\widehat{\grads u}\|_{L^2(\R^d)}\\
        &=
        \|(2\pi)^{r+s}|\cdot|^{r+s}\widehat{u}\|_{L^2(\R^d)}\\
        &=
        |u|_{H^{s+r}(\R^d)}.
    \end{aligned}
    \]
On the other hand,
\[ \begin{split}
\| \grads u \|_{L^2(\R^d)}^2 & = 2 \pi \int_{\R^d} |\xi|^{2s} |\widehat{u}(\xi)|^{2} \, d\xi \\
& \le 2 \pi \left( \int_{B_1} |\widehat{u}(\xi)|^{2} \, d\xi + 
\int_{B_1^c} |\xi|^{2(s+r)} |\widehat{u}(\xi)|^{2} \, d\xi \right) \\
& \le 2\pi \left( \| u \|_{L^2(\R^d)}^2 + | u |_{H^{s+r}(\R^d)}^2 \right).
\end{split} \]
\end{proof}

We summarize the preceding discussion about regularity of solutions in the following proposition.
\begin{proposition}[regularity of solutions]\label{prop:regularity-solutions}
 Let $\Omega \subset \R^d$ be a bounded, Lipschitz domain and $f \in L^2(\Omega)$. Consider $ (p, \vPhi)$ the weak solution of \eqref{eq:Darcy}. We have
\begin{equation}\label{eq:regularity_s_greater_12}
    \begin{aligned}
      \|p\|_{{\widetilde H}^{s+\frac12-\eps}(\Omega)}+|\vPhi|_{H^{\frac12-\eps}(\R^d)}&\leq\frac{C}{\sqrt{\eps|1-2s|}}\|f\|_{L^2(\Omega)},\quad &\text{for $s>\frac12$},\\
        \|p\|_{{\widetilde H}^{2s-\eps}(\Omega)}+|\vPhi|_{H^{s-\eps}(\R^d)}&\leq\frac{C}{\sqrt{\eps|1-2s|}}\|f\|_{L^2(\Omega)},\quad &\text{for $s<\frac12$},\\
        \|p\|_{{\widetilde H}^{1-\eps}(\Omega)}+|\vPhi|_{H^{\frac12-\eps}(\R^d)}&\leq \frac{C}{\eps}\|f\|_{L^2(\Omega)},\quad &\text{for $s=\frac12$},
    \end{aligned}
    \end{equation}
    for any $\eps<\min\{2s,\frac12+s\}$, with a constant $C = C(d,\Omega)$.
\end{proposition}

\begin{remark}[regularity in case $f \in B^{-s+1/2}_{2,1}(\Omega)$]
We stated the previous proposition under the simplifying assumption $f \in L^2(\Omega)$. However, as we already commented, the technique from \cite{Borthagaray2023} is better suited for $f$ belonging  to a certain Besov space. It turns out that, if $s > 1/2$, the assumption $f \in L^2(\Omega)$ in Proposition \ref{prop:regularity-solutions} can be relaxed to $f \in B^{-s+1/2}_{2,1}(\Omega)$. Additionally, if $s \le 1/2$, then assuming $f \in B^{-s+1/2}_{2,1}(\Omega)$ is stronger than assuming $f \in L^2(\Omega)$, but one also has the stronger estimate
\begin{equation} \label{eq:max-regularity}
   \|p\|_{{\widetilde H}^{s+\frac12-\eps}(\Omega)}+|\vPhi|_{H^{\frac12-\eps}(\R^d)}\leq\frac{C}{\sqrt{\eps}}\|f\|_{B^{-s+1/2}_{2,1}(\Omega)}.
\end{equation}

See \cite[Corollary 2.1]{Barrett}; we also refer to that work for a short discussion on Besov spaces.
\end{remark}

The inequality \eqref{eq:max-regularity} highlights that Besov regularity of the right-hand side translates into improved Sobolev-type estimates for the solution. In what follows, we pursue an analogous result involving weighted Bessel-type norms, where the regularity is expressed in terms of the Hölder continuity of the data.

Indeed, Hölder  regularity of solutions to \eqref{eq:Poisson} is well-known for Lipschitz domains and bounded right-hand sides, and was established in \cite{ros2014dirichlet}. Furthermore, when $f$ has some additional Hölder regularity one can deduce estimates for $u$ in certain weighted Hölder norms. 
Let $\delta(x)=d(x,\partial\Omega)$ and $\delta(x,y)=\min\{\delta(x),\delta(y)\}$.
For $\beta>0$ with $\beta=k+\beta'$, $k\in\N$ and $\beta'\in[0,1)$, and $\theta\ge-\beta$ we define
\[
|u|_{C^{k,\beta'}(\Omega)}:=\sup_{x,y\in\Omega}\frac{|u(x)-u(y)|}{|x-y|^{\beta'}},
\]
i.e, the semi-norm of the Hölder  space $C^{k,\beta'}(\Omega)$. We also define
\[
|u|_{\beta}^{(\theta)}:=\sup_{x,y\in\Omega}\delta(x,y)^{\beta+\theta}\frac{|\nabla^k u(x)-\nabla^k u(y)|}{|x-y|^{\beta'}},
\]
together with the norm:
\begin{itemize}
    \item for $\theta\ge0$,
    \[
        \|u\|_{\beta}^{(-\theta)}:=\sum_{j=0}^k\|\delta^{j+\theta}\nabla^ju\|_{L^\infty(\Omega)}+|u|_{\beta}^{(\theta)},
    \]
    \item for $0>\theta\ge-\beta$,
    \[
        \|u\|_{\beta}^{(-\theta)}:=\|u\|_{C^{\lfloor -\theta\rfloor,-\theta\,+ \lfloor \theta\rfloor}(\Omega)}\sum_{j=1}^k\|\delta^{j+\theta}\nabla^ju\|_{L^\infty(\Omega)}+|u|_{\beta}^{(\theta)}.
    \]
\end{itemize}

With these definitions in place, we recall the following results on the Hölder regularity of solutions to \eqref{eq:Poisson}; see \cite[Propositions 1.1 and 1.4]{ros2014dirichlet}.

\begin{proposition}[weighted H\"older regularity]
    Let $\Omega\subset\R^d$ a Lipschitz domain satisfying the exterior ball condition and $f\in L^\infty(\Omega)$. Then the solution $u$ of \eqref{eq:Poisson} belongs to $C^{0,s}(\R^d)$ and 
    \begin{equation} \label{eq:cs-regularity}
    \|u\|_{C^{0,s}(\overline\Omega)}\leq C(\Omega,s)\|f\|_{L^\infty(\Omega)}.        
    \end{equation}
    Moreover, let
    $\beta>0$ such that neither $\beta$ nor $\beta+2s$ is an integers. Then, if $f\in C^{\beta}(\Omega)$ with $\|f\|_\beta^{(-s)} < \infty$, $u\in C^{\beta+2s}(\Omega)$ and
    \begin{equation} \label{eq:weighted-cs-regularity}
    \|u\|_{\beta+2s}^{(-s)}\leq C(\Omega,s,\beta)\left(\|u\|_{C^{0,s}(\Omega)}+\|f\|_\beta^{(-s)}\right).        
    \end{equation}
\end{proposition}

Let $\beta$ be such that
\begin{equation} \label{eq:cond-beta}
   \beta\in\begin{cases}
    (1,2-2s)\quad\text{if $s\in(0,\frac12)$,} \\
    (0,2-2s)\quad\text{if $s\in[\frac12,1)$}.
\end{cases} 
\end{equation}
Estimate \eqref{eq:weighted-cs-regularity} yields
\begin{equation}
\label{eq:ros1}
    |u|_{\beta+2s}^{(-s)}=\sup_{x,y\in\Omega}\delta(x,y)^{\beta+s}\frac{|\nabla u(x)-\nabla u(y)|}{|x-y|^{\beta+2s-1}}\leq C(\Omega,s,\beta,\|f\|_\beta^{(s)})
\end{equation}
and
\begin{equation}
\label{eq:ros2}
    \sup_{x\in\Omega}|\delta(x)^{1-s}\nabla u(x)|\leq C(\Omega,s,\beta,\|f\|_\beta^{(s)}).
\end{equation}

These two estimates allow us to bound $|\delta(x)^\alpha\grad^\gamma u(x)|$ for $x\in\Omega$. 

\begin{theorem}[weighted pressure estimates]
\label{thm:regularidad_con_pesos}
    Let $\Omega\subset\R^d$ a bounded Lipschitz domain satisfying the exterior ball condition, the function 
    $f\in C^{\beta}(\Omega)$ with $\|f\|_\beta^{(-s)} < \infty$
    for some $\beta \in [1-s, 2-2s) \setminus \{1\}$,
    and $(p,\vPhi)$ be the weak solution of \eqref{eq:Darcy}. Then, $p$ satisfies
    \[
        \|\delta^{\frac12-\eps}\grad^{1+s-2\eps}p\|_{L^2(\Omega)}\leq \frac{C(\Omega,s,f)}{(\beta-1+s+2\eps)(s-2\eps)\sqrt{\eps}}<\infty
    \quad
    \mbox{for all } \eps \in (0,s/2). 
     \]
\end{theorem}
\begin{proof}
 Let $\gamma\in(1,2)$, $\alpha>0$ and $\beta>0$ be such that 
 $f\in C^{\beta}(\Omega)$ with $\|f\|_\beta^{(-s)} < \infty$.
 Specific requirements on these parameters are derived in the following. According to \eqref{eq:def-grads}, for $x\in\Omega$ we have
    \[
    \begin{aligned}
    |\delta(x)^\alpha\grad^\gamma p(x)|
    &\leq 
C\int_{\R^d}\delta(x)^\alpha\frac{|p(x+h)-p(x)-\nabla p(x)\cdot h|}{|h|^{d+\gamma}}dh.
    \end{aligned}
    \]
    We split the integration domain as $\R^d=B(0,\delta(x)/2)\,\cup\,B(0,\delta(x)/2)^c$ and denote the corresponding integrals as $(I)$ and $(II)$, respectively.

    For $h\in B(0,\delta(x)/2)$,  by the mean value theorem we can write $p(x+h)-p(x)=\nabla p(x+\lambda(h)\,h)\cdot h$ for some $\lambda\in(0,1)$. Moreover, we have $\delta(x)\simeq\delta(x,x+\lambda(h)\,h)$. Therefore, if $\beta+2s-\gamma>0$, estimate \eqref{eq:ros1} gives
    \[
    \begin{aligned}
    (I)&\le C \delta(x)^\alpha\int_{B(0,\delta(x)/2)}\frac{|\nabla p(x+\lambda(h)\,h)-\nabla p(x)|}{|h|^{d+\gamma-1}}dh\\ 
    &\leq
    C(\Omega,s,\beta,f)\,\delta(x)^{\alpha-\beta-s}\int_{B(0,\delta(x)/2)}\frac{dh}{|h|^{d+\gamma-\beta-2s}}\\
    &\leq
    \frac{C(\Omega,s,\beta,f)}{\beta+2s-\gamma}\delta(x)^{\alpha-\gamma+s}.
    \end{aligned}
    \]

    In order to bound $(II)$ we employ the Hölder regularity  \eqref{eq:cs-regularity} and the estimate \eqref{eq:ros2} to deduce 
    \[
    \begin{aligned}
        (II)&\leq C\left[\delta^\alpha(x)\int_{B(0,\delta(x)/2)^c}\frac{1}{|h|^{d+\gamma-s}}dh+\delta(x)^{\alpha-1+s}\int_{B(0,\delta(x)/2)^c}\frac{1}{|h|^{d+\gamma-1}}dh\right]\\
        &\leq
        \frac{C(\Omega,s,\beta,f)}{\gamma-1}\delta(x)^{\alpha-\gamma+s},
    \end{aligned}
    \]
    because we are assuming $\gamma>1$. Therefore, we conclude
    \begin{equation}
    \label{eq:weighted-grad-pointwbound}
    |\delta(x)^\alpha\grad^\gamma p(x)|\leq \frac{C(\Omega,s,\beta,f)}{(\beta+2s-\gamma)(\gamma-1)}\delta(x)^{\alpha-\gamma+s}
    \end{equation}
    for $\beta+2s-\gamma>0$.  Taking squares and integrating over $\Omega$, we obtain
    \[
    \begin{aligned}
    \|\delta^\alpha\grad^\gamma p\|_{L^2(\Omega)}^2&\leq\frac{C(\Omega,s,\beta,f)}{(\beta+2s-\gamma)^2(\gamma-1)^2}\int_\Omega\delta(x)^{2(\alpha-\gamma+s)}dx\\
    &\leq
    \frac{C(\Omega,s,\beta,f)}{(\beta+2s-\gamma)^2(\gamma-1)^2(1-2(\alpha+s-\gamma))},
    \end{aligned}
    \]
    where we have used \cite[Lemma 2.14]{Mazya} and assumed 
    \begin{equation} \label{eq:cond-parameters}
            \beta>\gamma-2s,\quad \alpha>-\frac12-s+\gamma.
    \end{equation}
As long as these two conditions are met, we can guarantee $\delta^\alpha \grad^\gamma p \in L^2(\Omega)$. 
    We choose, for $\eps>0$, $\alpha=\frac12-\eps$ and $\gamma=1+s-2\eps$ to conclude
    \[
        \|\delta^{\frac12-\eps}\grad^{1+s-2\eps} p\|_{L^2(\Omega)}\leq \frac{C(\Omega,s,f)}{(\beta-1+s+2\eps)(s-2\eps)\sqrt{\eps}}.
    \]
\end{proof}

\begin{remark}[global estimates]
The same technique allows one to show that, for $\Omega_\lambda$ being a neighborhood of $\Omega$, one has $\delta^{\frac12-\eps}\grad^{1+s-2\eps} p \in L^2(\Omega_\lambda)$. Indeed, it suffices to observe that, for $x \in \Omega_\lambda \setminus \Omega$, it holds that
\[
|\grad^{1+s-2\eps} p (x)| \le C \int_{B(0,\delta(x))^c} \frac{|p(x+h)|}{|h|^{d+1+s-2\eps}} \, dh ,
\]
and therefore the same argument as for the term (II) in the previous proof can be applied. On regions uniformly away from $\partial \Omega$, this bound also shows that $\grad^{1+s-2\eps}$ is smooth; see Lemma \ref{lemma:flux-interpolation} below.
\end{remark}

The regularity provided by Theorem \ref{thm:regularidad_con_pesos} is analogous to the fractional weighted Sobolev regularity established in \cite[Proposition 3.12]{acosta2017fractional}. 
This suggests that graded meshes could be employed to exploit such regularity and achieve higher convergence rates in the numerical scheme. However, proving that these improved rates are indeed attained would require a weighted (or enhanced) Poincaré inequality in the spirit of \cite[Proposition 4.7]{acosta2017fractional}, which remains a subject of ongoing research. 
Nevertheless, our numerical experiments in Section 6 exhibit improved convergence orders when graded meshes are used.

\begin{remark}[choice of parameters]
The choice of parameters in Theorem \ref{thm:regularidad_con_pesos} is motivated by its application for finite element approximations in $d=2$ dimensions. Indeed, one can show $\delta^\alpha \grad^\gamma p \in L^2(\Omega)$ as long as \eqref{eq:cond-parameters} is satisfied, but the application of weighted regularity estimates requires the use of adapted meshes. For the discretization of the Poisson problem \eqref{eq:Poisson} in $d=2$ dimensions, reference \cite{BoCi18} shows that the choice of parameters we made in Theorem \ref{thm:regularidad_con_pesos} gives rise to optimal interpolation error bounds in $\widetilde{H}^s(\Omega)$, with respect to the number of degrees of freedom, over shape-regular meshes.
\end{remark}

\section{Finite element discretization} \label{sec:FE}

By replacing $\calL$ with $\calL_{\rm stab}$, we obtain a coercive formulation. Consequently, any conforming finite element space yields a stable discretization. In this work, we focus on linear Lagrange elements for both the pressure $p$ and the flux $\vPhi$.

Let us begin by describing the discrete framework that we will use. We observe that we are approximating $\vPhi$, which is not compactly supported, and both the form $a$ in \eqref{eq:def-ab} and the stabilization term in \eqref{eq:def-calLstab2} involve integration in $\R^d$. To tackle this issue, our computational domain will be a ball $B_H$ with radius $H = H(h)\geq1$, containing $\Omega$ and such that $H \simeq d(\overline{\Omega},B_H^c)$.

Let $\left\{ \Th(B_H) \right\}_{h>0}$ be a family of simplicial meshes of $\overline{B_H}$, whose elements $\{T \}_{T \in \Th(B_H)}$ are assumed to be closed, that satisfy:
\begin{itemize}
    \item Shape-regularity, i.e., there exist a constant $\sigma>0$ such that 
    \begin{equation}
    \label{eq:shape-regularity}
    \sup_{h>0}\sup_{T\in\Th}\frac{h_T}{\rho_T}=\sigma,
    \end{equation}
    where $h_T=diam(T)$ and $\rho_T$ is the diameter of the largest ball contained in $T$.
    \item For every $h>0$, the set 
    \[\Th(\Omega)=\{T\in\Th(B_H):T\cap\Omega\not=\emptyset\}\]
    is a simplicial triangulation of $\overline{\Omega}$.
\end{itemize}
We denote  by $\calN_h = \{ z_i \colon i = 1, \ldots , N_h \}$ the set of nodes of $\Th(B_H)$. We set $n_h=\#(\calN_h\cap\Omega)$ and assume that the nodes are labeled so that those belonging to $\Omega$ come first, i.e., $\calN_h\cap\Omega = \{z_1 , \ldots, z_{n_h}\}$.
Let $\{\varphi_i\}_{i=1}^{N_h}$ be the standard piecewise linear Lagrange nodal basis associated with $\calN_h$ and $B_i$ be the largest ball centered at $z_i$ and contained in $supp(\varphi_i)$. The finite element space is defined as
\begin{equation}\label{eq:def-FE-space}
    \bV_h = \{ (q_h,\vPsi_h)\in\calP_1(\Th(B_H))\times[\calP_1(\Th(B_H))]^d\subset\bV \, \colon q_h|_{\Omega^c} = 0 \},
\end{equation}
where we assume that discrete functions are extended by zero outside of the computational domain $B_H$.

The discretization of problem \eqref{eq:Darcy-stabilized} reads: find $(p_h, \vPhi_h) \in \bV_h$ such that
\begin{equation} \label{eq:Darcy-stabilized-weak}
 \calL_{\rm stab} ((p_h, \vPhi_h), (q_h, \vPsi_h)) = F(q_h) \quad  \forall  (q_h, \vPsi_h) \in \bV_h.
\end{equation}

The coercive formulation and the fact that $\bV_h\subset\bV$ immediately imply the existence and uniqueness of solutions to \eqref{eq:Darcy-stabilized-weak}, and the best approximation property stated below.
\begin{proposition}[best approximation] \label{prop:best-approximation}
Let $(p, \vPhi) \in \bV$ and $(p_h, \vPhi_h) \in \bV_h$ be the solutions to \eqref{eq:Darcy-stabilized} and \eqref{eq:Darcy-stabilized-weak}, respectively.
We have the following Galerkin orthogonality:
\begin{equation}\label{eq:orthogonality}
 \calL_{\rm stab} ((p-p_h, \vPhi-\vPhi_h), (q_h, \vPsi_h)) = 0 \quad  \forall  (q_h, \vPsi_h) \in \bV_h.
\end{equation}
Consequently, we obtain
\begin{equation}\label{eq:best-approximation}
\vertiii{(p-p_h, \vPhi-\vPhi_h)} = \min_{(q_h, \vPsi_h) \in \bV_h} \vertiii{(p-q_h, \vPhi-\vPsi_h)}.
 \end{equation}
\end{proposition}
\begin{proof}
    Let $(q_h, \vPsi_h) \in \bV_h\subset\bV$. First, we have
    \begin{align*}
        \calL_{\rm stab} ((p-p_h, \vPhi-\vPhi_h), (q_h, \vPsi_h)) &= \calL_{\rm stab} ((p, \vPhi), (q_h, \vPsi_h))-\calL_{\rm stab} ((p_h,\vPhi_h), (q_h, \vPsi_h))\\
        &=F(q_h)-F(q_h)\\
        &=0.
    \end{align*}
    In second place, by the coercivity and continuity of $\calL_{\rm stab}$ and the Galerkin orthogonality, we deduce that, for all $(q_h, \vPsi_h) \in \bV_h$,
    \begin{align*}
        \vertiii{(p-p_h,\vPhi-\vPhi_h)}^2&=\calL_{\rm stab} ((p-p_h, \vPhi-\vPhi_h), (p-p_h, \vPhi-\vPhi_h))\\
        &=\calL_{\rm stab} ((p-p_h, \vPhi-\vPhi_h), (p-q_h, \vPhi-\vPsi_h))\\
        &\le\vertiii{(p-p_h,\vPhi-\vPhi_h)} \, \vertiii{(p-q_h,\vPhi-\vPsi_h)},
    \end{align*}
and \eqref{eq:best-approximation} follows.
\end{proof}

\subsection{Quasi-interpolation}
Our next task to derive convergence rates is to obtain interpolation estimates.
The standard Lagrange interpolation is not a feasible option in our setting because of the low regularity of solutions and its lack of stability in the corresponding low-order fractional Sobolev spaces. Therefore, we will use the quasi-interpolation operator introduced in \cite{chen2000residual}. Other suitable choices of quasi-interpolation (e.g. Cl\'ement, Scott-Zhang) would also be adequate for our purposes.

\begin{definition}[quasi-interpolation operators]
We define $\inte:L^1(\Omega) \to\calP_1(\Th(B_H))$ and $\iinte:(L^1(\Omega))^d\to\calP^d_1(\Th(B_H))$ as 
\begin{equation}
\begin{aligned}
&\inte q:=\sum_{z_i\in\Omega}\left(\frac{1}{|B_i|}\int_{B_i}q(x)\,dx\right)\varphi_i,\\
&\iinte\vPsi:=\sum_{z_i\in B_H}\left(\frac{1}{|B_i|}\int_{B_i}\vPsi(x)\,dx\right)\varphi_i.
\end{aligned}
\end{equation}
\end{definition}
We refer to \cite{chen2000residual,borthagaray2019weightedsobolevregularityrate} for basic properties of this operator. 
We are concerned with its stability and approximation properties with respect to fractional-order seminorms.
To this end, given $T\in\Th(B_H)$, we define the sets
\[
    \begin{aligned}
        & S^1_T = \bigcup_{T'\in \Th(B_H):\,T'\,\cap\, T\not=\emptyset} T', \quad
        & S^2_T = \bigcup_{T'\in \Th(B_H):\,T'\,\cap\, S^1_T\not=\emptyset} T'.
    \end{aligned}
\]
It is well known that fractional semi-norms are, in general, not additive with respect to domain decompositions. However, by allowing some overlapping, localization becomes possible. In \cite{faermann2000localization}, the following localization of the Gagliardo semi-norm is provided:
\begin{equation}
\label{eq:localization}
|q|^2_{H^s(\Omega)}\leq C\sum_{T\in\Th(\Omega)}\left(\iint_{T\times (S^1_T\cap\Omega)}\frac{|q(x)-q(y)|^2}{|x-y|^{d+2s}}dy\,dx + \frac{\|q\|_{L^2(T)}^2}{h_T^{2s}}\right).
\end{equation}
The constant above depends on $d, s$, and the shape-regularity parameter $\sigma$ from \eqref{eq:shape-regularity}.
We are interested in an estimate of this kind for the $\tHs$-norm, i.e, for $|q|_{H^s(\R^d)}$ when $q\in\tHs$. This is the goal of the next lemma (see also \cite[Lemma 4.1]{borthagaray2023constructive}).

\begin{lemma}[localization of the $\tHs$-norm]
\label{lema:localization-in-tHs}
Let $\Omega$ be a bounded Lipschitz domain and $\Th(B_H)$ as above. Then it holds,
\[
\|q\|_{\tHs}^2 = |q|_{H^s(\R^d)}^2 \leq C \sum_{T\in\Th(\Omega)}\left(\iint_{T\times S^1_T}\frac{|q(x)-q(y)|^2}{|x-y|^{d+2s}}dy\,dx + \frac{\|q\|_{L^2(T)}^2}{h_T^{2s}}\right),
\]
for all $q\in\tHs$.
\end{lemma}
\begin{proof}
    Let $\Omega_h=\bigcup_{T\in\Th(\Omega)}S^1_T$ and $q\in\tHs$. Observing that the integral over $\Omega_h^c\times\Omega_h^c$ in the definition of $|q|_{H^s(\R^d)}$ is zero, we obtain
    \begin{equation}
    \label{eq:prueba_localizacion2}
    \begin{aligned} \frac{2}{\nu(d,s)}
    |q|_{H^s(\R^d)}^2
    &=
    \iint_{\Omega_h\times\Omega_h}\frac{|q(x)-q(y)|^2}{|x-y|^{d+2s}}dy\,dx+2\iint_{\Omega_h\times\Omega_h^c}\frac{|q(x)-q(y)|^2}{|x-y|^{d+2s}}dy\,dx\\
    &\le
    \frac{2}{\nu(d,s)} |q|_{H^s(\Omega_h)}+2\int_{\Omega}|q(x)|^2\int_{\Omega_h^c}|x-y|^{-d-2s}dy\,dx.
    \end{aligned}
    \end{equation}
Now, for every element $T \subset \Omega$, by the shape-regularity of $\Th$ we have $d(T,\Omega_h^c) \ge C h$. Thus, integrating in polar coordinates we derive
\[
\int_{\Omega}|q(x)|^2\int_{\Omega_h^c}|x-y|^{-d-2s}dy\,dx \le C \sum_{T \in \Th(\Omega)} \frac{\|q\|_{L^2(T)}^2}{h_T^{2s}},
\]
that yields
\begin{equation}\label{eq:prueba_localizacion}
|q|_{H^s(\R^d)}^2 \le C \left(|q|^2_{H^s(\Omega_h)}+\sum_{T \in \Th(\Omega)} \frac{\|q\|_{L^2(T)}^2}{h_T^{2s}} \right).    
\end{equation}

    Next, using \eqref{eq:localization} we deduce
    \[
    |q|_{H^s(\Omega_h)}^2\leq C(d,s) \sum_{T\subset\Omega_h}\left(\iint_{T\times S^1_T}\frac{|q(x)-q(y)|^2}{|x-y|^{d+2s}}dy\,dx + \frac{\|q\|_{L^2(T)}^2}{h_T^{2s}}\right).
    \]
    The sum above contains many terms that are zero, as interactions between elements in $\Omega^c$ vanish. Indeed, $q\equiv0$ in $\Omega^c$, and hence for elements $T,T'\subset\Omega_h\backslash\Omega$ it holds
    \[
    \iint_{T\times T'}\frac{|q(x)-q(y)|^2}{|x-y|^{d+2s}}dy\,dx=\iint_{T'\times T}\frac{|q(x)-q(y)|^2}{|x-y|^{d+2s}}dy\,dx=0.
    \]
    Thus,
    \[
    \iint_{T\times S^1_T}\frac{|q(x)-q(y)|^2}{|x-y|^{d+2s}}dy\,dx=\iint_{T\times (S^1_T\,\cap\, \Omega)}\frac{|q(x)-q(y)|^2}{|x-y|^{d+2s}}dy\,dx,
    \]
    for every $T\subset\Omega_h$. Therefore, we deduce
    \[
    |q|_{H^s(\Omega_h)}\leq C(d,s)\sum_{T\in\Th(\Omega)}\left(\iint_{T\times S^1_T}\frac{|q(x)-q(y)|^2}{|x-y|^{d+2s}}dy\,dx + \frac{\|q\|_{L^2(T)}^2}{h_T^{2s}}\right).
    \]
The proof follows by combining this estimate with
 \eqref{eq:prueba_localizacion}.
\end{proof}

Lemma \ref{lema:localization-in-tHs} allows to obtain global interpolation estimates on $\tHs$ from local considerations.
Let $s\in(0,1)$, $t\in(s,2]$ and $\sigma$ the shape-regularity constant \eqref{eq:shape-regularity}. For $T\in\Th(B_H)$, we have
\begin{equation}
\label{eq:local-interpolation-estimates}
\iint_{T\times S^1_T}\frac{|(q-\inte q)(x)-(q-\inte q)(y)|^2}{|x-y|^{d+2s}}dy\,dx\leq \frac{C(d,\sigma,t)}{1-s}h_T^{2(t-s)}|q|^2_{H^t(S^2_T)},
\end{equation}
see \cite[Proposition 4.10]{borthagaray2019weightedsobolevregularityrate}. 
Additionally, for $t\in[0,2]$ we have the $L^2$-approximation bound
\begin{equation}
\label{eq:L2-stability-interpolator}
\|q-\inte q\|_{L^2(T)}\leq Ch_T^t|q|_{H^t(S^1_T)},
\end{equation}
with a constant $C$ independent of $T$ and $h$. For $t=0$, the inequality follows from \cite[Lemma 3.1]{chen2000residual}, while for $t=2$ the inequality follows from \cite[Lemma 3.2]{chen2000residual}. By interpolation between the cases $t=0$ and $t=2$, we obtain \eqref{eq:L2-stability-interpolator}. 
Naturally, estimates \eqref{eq:local-interpolation-estimates} and \eqref{eq:L2-stability-interpolator} also hold for $\iinte$.

Combining the local interpolation estimates \eqref{eq:local-interpolation-estimates} and \eqref{eq:L2-stability-interpolator}, and the localization provided by Lemma \ref{lema:localization-in-tHs}, we deduce global interpolation estimates. In particular, for quasi-uniform meshes, these read as follows.
\begin{lemma}[global interpolation estimates]
    Let $s\in(0,1)$, $t\in(s,2]$, $r \in (0,2]$, and $\Th(B_H)$ be a quasi-uniform mesh. Then, we have
    \begin{equation}\label{eq:interpolation-estimates}
        \begin{aligned}
            &\|q-\inte q\|_{\tHs}\leq C(d,\sigma,t)h^{t-s}\|q\|_{{\widetilde H}^{t}(\Omega)},\\
            &\|\vPsi-\iinte\vPsi\|_{L^2(B_H)}\leq C(d,\sigma,t)h^r|\vPsi|_{H^r(\R^d)},
        \end{aligned}
    \end{equation}
    for all $q\in{\widetilde H}^{t}(\Omega)$ and $\vPsi\in [H^r(\R^d)]^d$.
\end{lemma}

The coercivity norm $\vertiii{(p,\vPhi)}$ involves the $\tHs$-seminorm of $p$ and the $L^2(\R^d)$ norm of $\vPhi$. Thus, taking into account \eqref{eq:interpolation-estimates}, to conclude an interpolation estimate we need to address $\| \vPhi - \iinte\vPhi \|_{L^2(B_H^c)}$. Since $\iinte\vPhi$ vanishes outside such a ball, this calculation reduces to bounding the decay of $|\vPhi|$. For the sake of simplicity, from this point on we assume that $B_H$ is centered at the origin.

\begin{lemma}[flux decay]\label{lemma:flux-interpolation}
Let $\Omega \subset \R^d$ be a bounded, Lipschitz domain such that $0 \in \Omega$, and $(p, \vPhi) \in \bV$ be the solution to \eqref{eq:Darcy-stabilized}. We have, for all $x \in \Omega^c$ and all  multi-index $\alpha\in\N^d$,
\[
|\partial^\alpha\vPhi (x)| \le \frac{C(d,s,\alpha)}{d(x,\Omega)^{d+s+|\alpha|}} \| p \|_{L^1(\Omega)} .
\]
Consequently,
\begin{equation}
\label{eq:flux-l2-bound}
\begin{aligned}
    &\|\partial^\alpha\vPhi\|_{L^2(B_H^c)}\leq C(d,s,\alpha,\Omega)H^{-\frac d2 -s-|\alpha|}\,\|f\|_{L^2(\Omega)}.
\end{aligned}
\end{equation}
\end{lemma}
\begin{proof}

First, for $z\in\R^d\backslash\{0\}$ it holds  
\[
\frac{\partial^\alpha}{\partial z}\left(\frac{z}{|z|^{d+s+1}}\right)=\frac{Q_\alpha(z)}{|z|^{d+s+1+2|\alpha|}},
\]
with $Q_\alpha:\R^d\to\R^d$ a polynomial vector field that is componentwise homogeneous of degree $|\alpha|+1$. Therefore, for $c\in\R$ depending only on the coefficients of the components of $Q_\alpha$, we have
\[
\left|\frac{\partial^\alpha}{\partial z}\left(\frac{z}{|z|^{d+s+1}}\right)\right|\le \frac{c}{|z|^{d+s+|\alpha|}}.
\]

Since $\vPhi = \grads p$ in $\R^d$, for all $x \in \Omega^c$ and $y \in \Omega$, we exploit the fact that $|x-y| \ge d(x,\Omega)$ to deduce
\[
\begin{aligned}
|\partial^\alpha\vPhi(x)|^2
&=
\left|\mu(d,s)\int_\Omega p(y)\frac{\partial^\alpha}{\partial x}\left(\frac{x-y}{|x-y|^{d+s+1}}\right)dy\right|^2\\
&\le
C\left(\int_\Omega |p(y)|\frac{1}{|x-y|^{d+s+|\alpha|}}dy\right)^2
 \\ 
&\le
\frac{C}{d(x,\Omega)^{2(d+s+|\alpha|)}} \| p \|_{L^1(\Omega)}^2,
\end{aligned}
\]
with $C=C(d,s,\alpha)$. 

Therefore, noticing that $d(x,\Omega)\geq c|x|$ for some constant $c$, a change of variables to polar coordinates gives the bound
\begin{equation}
\begin{aligned}
\|\partial^\alpha\vPhi\|_{L^2(B^c_H)}^2&\leq C\,\|p\|_{L^1(\Omega)}^2\int_{B_H^c}\frac{1}{d(x,\Omega)^{2(d+s+|\alpha|)}}\,dx\\
&= 
C\,\|p\|_{L^1(\Omega)}^2\int_{H}^\infty \rho^{-d-2s-2|\alpha|-1}\,d\rho\\
&\leq
C\,\|p\|_{L^1(\Omega)}^2 H^{-d-2s-2|\alpha|},
\end{aligned}
\end{equation}
with a constant $C=C(d,s,\alpha)$. The inequality \eqref{eq:flux-l2-bound} follows by 
observing $\|p\|_{L^1(\Omega)}\leq C(\Omega,s) \|f\|_{L^2(\Omega)}$.
\end{proof}

\subsection{Convergence in the stability norm} 
Collecting the previous interpolation and decay estimates, together with the regularity of solutions, we obtain convergence rates for our numerical scheme.

\begin{proposition}[order of convergence]\label{cor:order-convergence}
Let $(p, \vPhi) \in \bV$ and $(p_h, \vPhi_h) \in \bV_h$ be the solutions to \eqref{eq:Darcy-stabilized} and \eqref{eq:Darcy-stabilized-weak}, respectively. Assume that $\Th(B_H)$ is a quasi uniform mesh and that the auxiliary mesh grows according to $H \geq (h |\log h|)^{\frac{-1}{d+2s}}$.

Then, we have the following convergence rates with respect to the mesh size $h>0$,
\begin{equation}\label{eq:order-convergence}
\vertiii{(p-p_h,\vPhi-\vPhi_h)}\leq\begin{cases} 
     Ch^\frac12|\log h|^\frac12\|f\|_{L^2{(\Omega)}},\quad&\text{for $s>\frac12$},\\
     Ch^\frac12|\log h|\,\|f\|_{L^2{(\Omega)}},\quad&\text{for $s=\frac12$},\\
     Ch^s|\log h|^\frac12\|f\|_{L^2{(\Omega)}},\quad&\text{for $s<\frac12$},
\end{cases}
\end{equation}
with a constant $C=C(d,s,\sigma)$. In particular, in terms of the number of interior nodes $n_h\in\N$, the previous order of convergence reads
\begin{equation}
    \vertiii{(p-p_h,\vPhi-\vPhi_h)}\leq\begin{cases} 
     Cn_h^{-\frac{1}{2d}}|\log n_h|^\frac12\|f\|_{L^2{(\Omega)}},&\quad\text{for $s>\frac12$},\\
     Cn_h^{-\frac{1}{2d}}|\log n_h|\,\|f\|_{L^2{(\Omega)}},&\quad\text{for $s=\frac12$},\\
     Cn_h^{-\frac{s}{d}}|\log n_h|^\frac12\|f\|_{L^2{(\Omega)}},&\quad\text{for $s<\frac12$}.
\end{cases}
\end{equation}
\end{proposition}
\begin{proof}
We provide details for the case $s>\frac12$. The other ones follow by the same argument. 
    By combining the best approximation property (cf. Proposition \ref{prop:best-approximation}) with the interpolation estimates \eqref{eq:interpolation-estimates} and Lemma \ref{lemma:flux-interpolation} with $\alpha = 0$, we obtain
    \begin{equation}
    \label{eq:prueba-ordenes}
    \vertiii{(p-p_h,\vPhi-\vPhi_h)}^2\leq C\left(h^{2t-2s}\|p\|_{{\widetilde H}^{t}(\Omega)}^2+h^{2r}|\vPhi|_{H^r(\R^d)}^2+H^{-d-2s}\|f\|_{L^2(\Omega)}^2\right),
    \end{equation}
    for $t,r\in(0,2)$.
    Now, by choosing $t=s+\frac12-\eps$ and $\eps = |\log h|^{-1}$ (so that $h^{-\eps}$ is constant for $h<1$), we can estimate the first term in the right hand side above by means of the regularity estimates (cf. Proposition \ref{prop:regularity-solutions});
    \begin{equation}
        h^{2t-2s}\|p\|_{{\widetilde H}^{t}(\Omega)}^2\leq Ch|\log h|\|f\|_{L^2{(\Omega)}}^2.
    \end{equation}
    The same argument with $r=\frac12-\eps$ and $\eps = |\log h|^{-1}$ shows
    \begin{equation}
        h^{2r}|\vPhi|_{H^r(\R^d)}^2\leq Ch|\log h|\|f\|_{L^2{(\Omega)}}^2.
    \end{equation}
    Finally, it suffices to observe that the third term in the right-hand side of \eqref{eq:prueba-ordenes} is (at least) 
    of the same order with respect to $h$ because of our choice $H$.

    The second estimate follows by observing that, for a quasi-uniform mesh, the number of interior nodes satisfies $n_h\simeq h^{-d}$.
\end{proof}

\begin{remark}[orders in terms of Besov regularity]
   For $f \in B^{-s+1/2}_{2,1}(\Omega)$, the proof of Proposition \ref{cor:order-convergence} can be repeated with the regularity estimate \eqref{eq:max-regularity} in place of \eqref{eq:regularity_s_greater_12}. In this case, e.g. the bound \eqref{eq:order-convergence} takes the form
   \begin{equation}
    \label{eq:order-convergence-besov}
    \vertiii{(p-p_h,\vPhi-\vPhi_h)}\leq Ch^\frac12|\log h|^\frac12\|f\|_{B^{-s+1/2}_{2,1}(\Omega)},
   \end{equation}
   for all $s\in(0,1)$.
\end{remark}

\begin{remark}
A straightforward calculation allowed us to use equation \eqref{eq:order-convergence} to obtain convergence rates in terms of the number of interior nodes. 
However, $n_h$ does not accurately reflect the size of the discrete problem. In fact, since we approximate $\R^d$ by meshing $B_H$ and assume $H \geq (h |\log h|)^{\frac{-1}{d+2s}}$, the use of uniform meshes implies that $N_h \gg n_h$, where $N_h$ denotes the total number of nodes in the mesh $\Th(B_H)$.
To write \eqref{eq:order-convergence} in terms of $N_h$, we observe that $|B_H|\simeq H^d$ and, for an quasi uniform mesh with size $h$, $|B_H|\simeq N_h h^{d}$, which implies 
$N_h\simeq \frac{H^d}{h^d}$,
and thus
\[
N_h^{-\frac{d+2s}{d(1+d+2s)}}\geq h.
\]
This, combined with the previous result, allows us to rewrite Proposition \ref{cor:order-convergence} in terms of the total degrees of freedom:
\begin{equation}\label{eq:order-convergence-dofs}
\vertiii{(p-p_h,\vPhi-\vPhi_h)}\leq\begin{cases} 
     CN_h^{-\frac{d+2s}{2d(1+d+2s)}}|\log N_h|^\frac12\|f\|_{L^2{(\Omega)}},\quad&\text{for $s>\frac12$},\\
     CN_h^{-\frac{d+1}{2d(2+d)}}|\log N_h|\,\|f\|_{L^2{(\Omega)}},\quad&\text{for $s=\frac12$},\\
     CN_h^{-\frac{s(d+2s)}{d(1+d+2s)}}|\log N_h|^\frac12\|f\|_{L^2{(\Omega)}},\quad&\text{for $s<\frac12$}.
\end{cases}
\end{equation}
\end{remark}

\subsection{Convergence order for the pressure in the $L^2$-norm}
According to the bound \eqref{eq:order-convergence}, the order of convergence for the pressure $p$ in \eqref{eq:Darcy} in the $\tHs$-norm is $\min\{\frac{1}{2},\,s\}$ (up to logarithmic factors) with respect to the mesh parameter $h$.
Here, we
develop an argument along the lines of the well-known Aubin-Nitsche duality trick to derive the order of convergence of $p$ in the $L^2$ norm. 

\begin{proposition}[Aubin-Nitsche argument] \label{prop:Aubin-Nitsche}
Let $(p,\vPhi) \in \bV$ be the solution to \eqref{eq:Darcy-stabilized} and $(p_h,\vPhi_h)\in\bV_h$ be the solution to \eqref{eq:Darcy-stabilized-weak}. Assume that $\Th(B_H)$ is a quasi uniform mesh and that the computational domain grows according to $H \geq (h |\log h|)^{\frac{-1}{d+2s}}$. Then, if $f \in B^{-s+1/2}_{2,1}(\Omega)$, it holds that
\begin{equation} \label{eq:error-bound-L2}
    \| p - p_h \|_{L^2(\Omega)} \le C h^{\frac{1}{2}+\min\{s,\frac12\}} | \log h|^{\kappa+\frac12} \| f \|_{B^{-s+\frac{1}{2}}_{2,1}(\Omega)},
\end{equation}
with a constant $C=C(d,s,\sigma)$. Here, $\kappa=1$ if $s=\tfrac12$ and $\kappa=\tfrac12$ otherwise.
\end{proposition}
\begin{proof}
Define the errors $e^p:=p-p_h$ and $e^\vPhi:=\vPhi-\vPhi_h$. From the Galerkin orthogonality \eqref{eq:orthogonality} we deduce 
\begin{equation}
\label{eq:ortogonalidad-aubin-nitsche}
    \int_{\R^d}\grads e^p\cdot(\grads q_h+\vPsi_h)=\int_{\R^d} e^\vPhi\cdot(\grads q_h-\vPsi_h),
\end{equation}
for all $(q_h,\vPsi_h)\in\bV_h$.
Let $(p^e,\vPhi^e) \in \bV$ be the weak solution to 
\begin{equation}
\label{eq:Poisson-Aubin-Nitsche}
\left\lbrace \begin{aligned}
 \vPsi + \grads q = 0 & \quad  \mbox{ in } \R^d, \\
\divs \vPsi = e^p & \quad  \mbox{ in } \Omega, \\
q = 0 &  \quad \mbox{ in } \Omega^c .
\end{aligned} \right.
\end{equation}
Then, using $\vPhi^e = -\grads p^e$, we obtain
\begin{align}
\label{eq:aubin-nitsche-problem}
    \|e^p\|_{L^2(\Omega)}^2 
    &= \int_{\R^d}\grads p^e \cdot \grads e^p\\
    &= \tfrac12\int_{\R^d}\grads p^e\cdot \grads e^p
      - \tfrac12\int_{\R^d}\vPhi^e\cdot \grads e^p .
\end{align}
Combining equations \eqref{eq:aubin-nitsche-problem} and \eqref{eq:ortogonalidad-aubin-nitsche}, with $q_h=\inte p^e$ and $\vPsi_h = -\iinte \vPhi^e$, we deduce
\[
\begin{split}
2\|e^p\|_{L^2(\Omega)}^2 = & \int_{\R^d} \grads e^p\cdot\grads [p^e-\inte p^e] - \int_{\R^d} \grads e^p\cdot(\vPhi^e-\iinte\vPhi^e)\\
 & + \int_{\R^d} e^\vPhi\cdot(\iinte\vPhi^e+\grads\inte p^e).\\
\end{split}
\] 
Since $\vPhi^e = -\grads p^e$, this simplifies to
\[ \begin{aligned}
2\|e^p\|_{L^2(\Omega)}^2
= & \int_{\R^d} \grads e^p 
   \cdot \big(\grads[p^e-\inte p^e] - (\vPhi^e-\iinte \vPhi^e)\big) \\
& - \int_{\R^d} e^\vPhi
   \cdot \big(\grads[p^e-\inte p^e] + (\vPhi^e-\iinte \vPhi^e)\big).
   \end{aligned}
\] 

Thus, we aim to bound the right hand side in the previous equality. First, observe that \eqref{eq:order-convergence-besov} implies
\[
\|\grads e^p\|_{L^2(\R^d)}+\|e^\vPhi\|_{L^2(\R^d)}\leq Ch^\frac12|\log h|^\frac12\|f\|_{B^{-s+1/2}_{2,1}(\Omega)}.
\]
Lastly, the same argument used in Proposition \ref{cor:order-convergence}, combining the interpolation estimates \eqref{eq:interpolation-estimates} and Lemma \ref{lemma:flux-interpolation} with $\alpha = 0$, but
applied to the dual problem \eqref{eq:Poisson-Aubin-Nitsche} gives the bound
\[
\|\grads[ p^e-\inte p^e]\|_{L^2(\R^d)}+\|\vPhi^e-\iinte\vPhi^e\|_{L^2(\R^d)}\leq Ch^{\min\{s,\frac12\}}|\log h|^{\kappa}\|e^p\|_{L^2(\Omega)},
\]
and \eqref{eq:error-bound-L2} follows.
\end{proof}

\subsection{On the meshing of $\Omega^c$}
We conclude this section by showing how to apply Lemma \ref{lemma:flux-interpolation} to construct meshes with less degrees of freedom than quasi uniform ones, while preserving the convergence rates \eqref{eq:order-convergence} and \eqref{eq:error-bound-L2}. Lemma \ref{lemma:flux-interpolation} implies that the flux $\vPhi$ has $H^2$ regularity in any element $T$ such that $d(T,\partial\Omega)>0$:
\[
|\vPhi|_{H^2(T)}\leq \frac{C(d,s)}{d(T,\Omega)^{d+s+2}}\|p\|_{L^1(\Omega)}h_T^\frac{d}{2}.
\]
We combine this with the $H^2$ interpolation estimate \eqref{eq:L2-stability-interpolator},
\[
\|\vPhi-\Pi_h\vPhi\|_{L^2(T)}\leq  C(d,\sigma)h_T^2|\vPhi|_{H^2(S^1_T)},
\]
to obtain
\[
\|\vPhi-\Pi_h\vPhi\|_{L^2(T)}\leq C(d,s,\sigma)\frac{h_T^{2+\frac{d}{2}}}{d(T,\Omega)^{d+s+2}}\|p\|_{L^1(\Omega)}.
\]
This inequality allows us to increase the diameter of elements sufficiently far from $\Omega$, thereby reducing the total number of degrees of freedom without compromising the global convergence rate. 
Given $h>0$, consider
\[\Omega_\alpha = \{x\in\R^d:d(x,\Omega)< h^\alpha\},\]
with $\alpha\in[0,1]$ to be determined.
We keep a quasi uniform mesh size $h$ within $\Omega_\alpha$ and aim to construct globally shape-regular meshes. This implies that meshes need to be locally quasi uniform and that we must avoid mismatches between the mesh sizes in $\Omega_\alpha$ and $\Omega_\alpha^c$, particularly for elements near $\partial\Omega_\alpha$. 

 In $\Omega_\alpha^c$, we set the diameter of the elements so that the global optimal order of Proposition \ref{cor:order-convergence} is preserved; for any $T\subset\Omega_\alpha^c$, it is sufficient to choose $h_T$ such that
\begin{equation}\label{eq:improved_meshes}
\frac{h_T^{2+\frac{d}{2}}}{d(T,\Omega)^{d+s+2}}\simeq h^\frac12.
\end{equation}

\begin{figure}[!htbp]
    \centering
\includegraphics[width=0.7\linewidth]{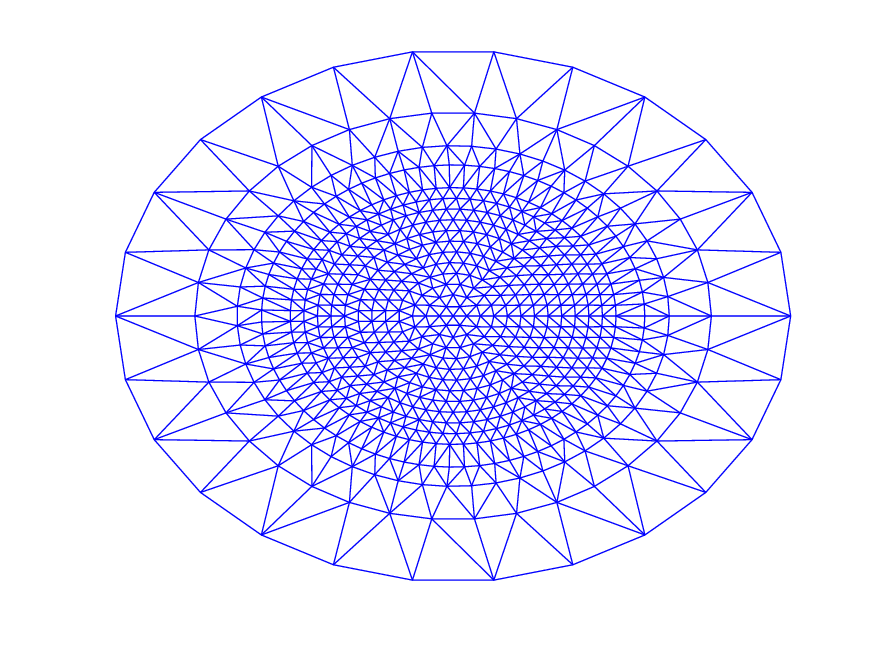}
    \caption{Example of a mesh satisfying \eqref{eq:construction_improved_meshes} for $\Omega=B(0,1)\subset\R^2$, $h=\frac1{10}$, and $\alpha = \frac59$, which corresponds to $s=\frac12$.} \label{fig:malla}
\end{figure}

Now, for elements $T\subset\Omega_\alpha^c$ with $T\cap\partial\Omega_\alpha^c\not=\emptyset$ we have $d(T, \partial\Omega)\simeq h^\alpha$. Substituting into \eqref{eq:improved_meshes}, we find
\[
h_T\simeq h^{\frac{1+2\alpha(d+s+2)}{4+d}}.
\]
Therefore, to ensure local quasi uniformity across the transition between $\Omega_\alpha$ and $\Omega_\alpha^c$ we need 
\[
\frac{1+2\alpha (d+s+2)}{4+d} = 1 \ \Rightarrow \ \alpha = \frac{3+d}{2(d+s+2)} .
\]
Summarizing, we construct the elements of $\Th(B_H)$ according to
\begin{equation}
\label{eq:construction_improved_meshes}
    h_T\simeq\begin{cases}
        h &\text{if $T\subset\Omega_\alpha$,}\\
        h^\frac{1}{4+d}d(T,\Omega)^{\frac{2(d+s+2)}{4+d}} & \text{if $T\subset\Omega_\alpha^c$,}
    \end{cases}
\end{equation}
with \[\Omega_\alpha = \{x\in\R^d:d(x,\Omega)< h^\alpha\} \quad \mbox{and} \quad 
\alpha = \frac{3+d}{2(d+s+2)}.\]
Figure \ref{fig:malla} displays a mesh constructed in this fashion for $\Omega = B(0,1) \subset \R^2$. 
Table \ref{tab:mesh_nodes} compares the total number of nodes required by a globally quasi-uniform mesh against a mesh satisfying \eqref{eq:construction_improved_meshes}. The results clearly demonstrate that the improved mesh strategy drastically reduces the number of degrees of freedom—by nearly one order of magnitude—while retaining the same convergence properties (see Table \ref{tab:errores_2d} below).

\begin{table}[!htbp]
    \centering
    \begin{tabular}{c|c|c|c|c} 
         & $h=0.1$ & $h=0.05$ & $h=0.025$ & $h=0.02$ \\
        \hline
Mesh satisfying \eqref{eq:construction_improved_meshes} & 701 & 2196 & 7378 & 11016 \\ \hline
Globally quasi-uniform mesh & 7541 & 29906 & 121214 & 187624 \\ 
\label{tab:mesh_nodes}
    \end{tabular}

    \caption{Comparison of the total number of nodes between a globally quasi-uniform mesh and a mesh satisfying \eqref{eq:construction_improved_meshes} with $s=\tfrac12$. Here, $\Omega = B(0,1)$, and the computational domain diameter is $H=2.72$.}
\end{table}

\section{Numerical experiments} \label{sec:experiments}
In this section, we present some experiments in dimensions $d=1$ and $d=2$.
The main challenges in computing solutions of \eqref{eq:Darcy} are the need to evaluate integrals over $\R^d$, the presence of a singular kernel, and the nonlocal nature of the problem, which leads to dense system matrices.  We begin this section with a brief discussion on the construction of the matrices involved in problem \eqref{eq:Darcy-stabilized-weak}, and refer to \cite{Acosta_2017} for further details on this aspect. 

We recall that, according to definition \eqref{eq:def-FE-space}, our discrete spaces consist of continuous functions $(q_h,\vPsi_h)\in\calP_1(\Th(B_H))\times[\calP_1(\Th(B_H))]^d\subset\bV$ such that $q_h$ vanishes on $\Omega^c$.
For the pressure, we consider the standard Lagrange nodal basis $B_{\text{pressure}}=\{\varphi_i\}_{1\leq i\leq n_h}$ defined on the interior nodes, while for the flux we employ the $d$-Lagrange nodal basis $B_{\text{flux}}=\{\vPhi_i\}_{1\leq i\leq dN_h}$. The difficult exercise when implementing \eqref{eq:Darcy-stabilized-weak} lies in the computation of the matrices $K \in \R^{n_h \times n_h}$, $B \in \R^{n_h\times d N_h}$,
\begin{equation}
    \begin{aligned}
    K_{ij}&=\int_{\R^d}\grads\varphi_i\cdot\grads\varphi_j \ dx=\frac{\mu(d,s)}{d+s-1}\iint_{\R^d\times\R^d}\frac{\nabla\varphi_i(y)\cdot\nabla\varphi_j(y)}{|x-y|^{d+s-1}}dy\ dx,\\
    B_{ij}&=\int_{\R^d}\grads\varphi_i\cdot\vPhi_j\,dx =\frac{\mu(d,s)}{d+s-1}\iint_{\R^d\times\R^d}\frac{\nabla\varphi_i(y)}{|x-y|^{d+s-1}}\cdot\vPhi_j(x)\ dy\ dx,
    \end{aligned}
\end{equation}
recall formula \eqref{eq:grads-divs-riesz-def}.
We point out here that, due to \eqref{eq:def-norma-Hs}, the matrix $K$ satisfies
\[
\begin{aligned}
K_{ij}
&=
\int_{\R^d}(-\Delta)^\frac s2\varphi_i(-\Delta)^\frac s2\varphi_j \ dx\\
&=\frac{\nu(d,s)}{2}\iint_{\R^d\times\R^d}\frac{(\varphi_i(x)-\varphi_i(y))(\varphi_j(x)-\varphi_j(y))}{|x-y|^{d+2s}}dx\,dy,
\end{aligned}
\]
so that it coincides with the stiffness matrix of problem \eqref{eq:Poisson}. Moreover, 
splitting 
$$\R^d\times\R^d=(B_H\times B_H)\cup(B_H^c\times B_H)\cup(B_H\times B_H^c)\cup(B_H^c\times B_H^c)$$ we obtain
\[ \begin{split}
K_{ij}  = & \frac{\nu(d,s)}{2}\iint_{B_H\times B_H}\frac{(\varphi_i(x)-\varphi_i(y))(\varphi_j(x)-\varphi_j(y))}{|x-y|^{d+2s}}dx\,dy\\ &+2\int_{B_H\times B_H^c}\frac{\varphi_i(x)\varphi_j(x)}{|x-y|^{d+2s}}dx\,dy.
\end{split}\]
In addition, observe that
\[
    B_{i,j}=\mu(d,s)\sum_{T_l\in sop\,\vPhi_j}\sum_{T_m\in sop\,\varphi_i}\nabla\varphi_i|_{T_m}\cdot\iint_{T_l\times T_m}\frac{\vPhi_j(x) }{|x-y|^{d+s-1}}dy\,dx.
\]
Hence, we iterate over the elements of the mesh using a double loop. More precisely, for any $T_l, T_m\in\Th(B_H)$ (Possibly with $T_l=T_m$), we must compute
\[
\nabla\varphi_i|_{T_m}\cdot\iint_{T_l\times T_m}\frac{\vPhi_j(x) }{|x-y|^{d+s-1}}dy\,dx,
\]
for the entries of $B$, and
\[
\iint_{T_l\times T_m}\frac{(\varphi_i(x)-\varphi_i(y))(\varphi_j(x)-\varphi_j(y))}{|x-y|^{d+2s}}dx\,dy,\quad
\int_{T_l\times B_H^c}\frac{\varphi_i(x)\varphi_j(x)}{|x-y|^{d+2s}}dx\,dy,
\]
for the entries of $K$.
The one-dimensional case is relatively simple, as the entries of both  matrices $K$ and $B$ can be explicitly computed even on non-uniform meshes.
The computation of $K$ in the two-dimensional case is thoroughly discussed in \cite{Acosta_2017} and we employ the same ideas to compute $B\in \R^{n_h\times 2N_h}$. The evaluation of the required integrals is carried out by mapping each pair of elements to a reference configuration and applying suitable quadrature rules. In particular, the quadrature rules described in \cite[Section 5.2]{sauter2011boundary} have the key advantage of transforming the integral over the product of two elements into an integral over $[0,1]^4$, where the singularities of the kernel can be explicitly computed. The extension of these techniques to the three-dimensional setting, for the matrix $K$, has been developed in \cite{MR4609814}.

The strategies described above allow us to assemble the system matrices required for the discretization of problem \eqref{eq:Darcy-stabilized-weak}. In the following, we test convergence rates of the stabilized method over quasi-uniform and graded meshes, as well as the influence of the computational domain diameter $H$ on the errors. To this end, it is convenient to consider test cases where exact solutions are available. In particular, explicit solutions of \eqref{eq:Poisson} are available when $\Omega$ is a ball. Consider the polynomials $P_{a,b,n}$ of degree $n$ defined as
\begin{equation}  \label{eq:Jacobi-pols}  
    P_{a,b,n}(t)=\frac{\Gamma(a+n+1)}{n!\,\Gamma(a+b+n+1)}\sum_{j=0}^n\binom{n}{j}\frac{\Gamma(a+b+n+j+1)}{\Gamma(a+j+1)}\left(\frac{t-1}{2}\right)^j,
\end{equation}
and the function $\varphi_s:\R^d\to\R$ such that
\[
    \varphi_s(x) = (1-|x|^2)^s\chi_{B(0,1)}.
\]
A family of solutions can be built by using this functions \cite[Theorem 3]{dyda2015fractionallaplaceoperatormeijer}.
\begin{theorem}
\label{thm:explicit_solutions}
    Let $\Omega=B(0,1) \subset \R^d$. For $s\in(0,1)$ and $n\in\N$, consider
    \[
        C(n,d,s) = \frac {n!\,\Gamma(\frac{d}{2}+n)} {2^{2s}\Gamma(1+s+n)\,\Gamma(\frac d2+s+n)} 
    \]
    and $p_{n,s}:\R^d\to\R$ such that
    \[
        p_{n,s}(x) = P_{s,\frac d2-1,n}(2|x|^2-1).
    \]
    Then, $u=C(n,d,s)  \varphi_s\,p_{n,s}$ is the solution of \eqref{eq:Poisson} with right-hand side $f=\,p_{n,s}$.
\end{theorem}
For $n=0$, Theorem \ref{thm:explicit_solutions} gives us an explicit solution for the fractional torsion problem,
\begin{equation} 
\label{eq:exp1}
\left\lbrace \begin{aligned}
(-\Delta)^s u = 1 & \quad  \mbox{ in } B(0,1), \\
u = 0 &  \quad \mbox{ in } B(0,1)^c,
\end{aligned} \right.
\end{equation}
which, in turn, corresponds to an explicit expression for $p$ in the mixed formulation: 
\begin{equation} \label{eq:Darcy-torsion}
\left\lbrace \begin{aligned}
 \vPhi + \grads p = 0 & \quad  \mbox{ in } \R^d, \\
\divs \vPhi = 1 & \quad  \mbox{ in } B(0,1), \\
p = 0 &  \quad \mbox{ in } B(0,1)^c .
\end{aligned} \right.
\end{equation}

\subsection{Quasi-uniform meshes}
As a first example, we analyze the convergence rates for this problem on quasi-uniform meshes in $d=1$ and $d=2$ dimensions. The parameter $H$ is chosen according to Proposition \ref{cor:order-convergence}, i.e, $H\geq (h|\log h|)^{\frac{-1}{d+2s}}$. For $d=1$, we consider a uniform mesh $\Th([-H,H])$ of mesh size $h$, and for $d=2$, we construct $\Th(B_H)$ following \eqref{eq:construction_improved_meshes}.

\begin{figure}[!htbp]
    \centering
    \begin{minipage}{0.49\textwidth}
    \centering
    \includegraphics[width=\textwidth]{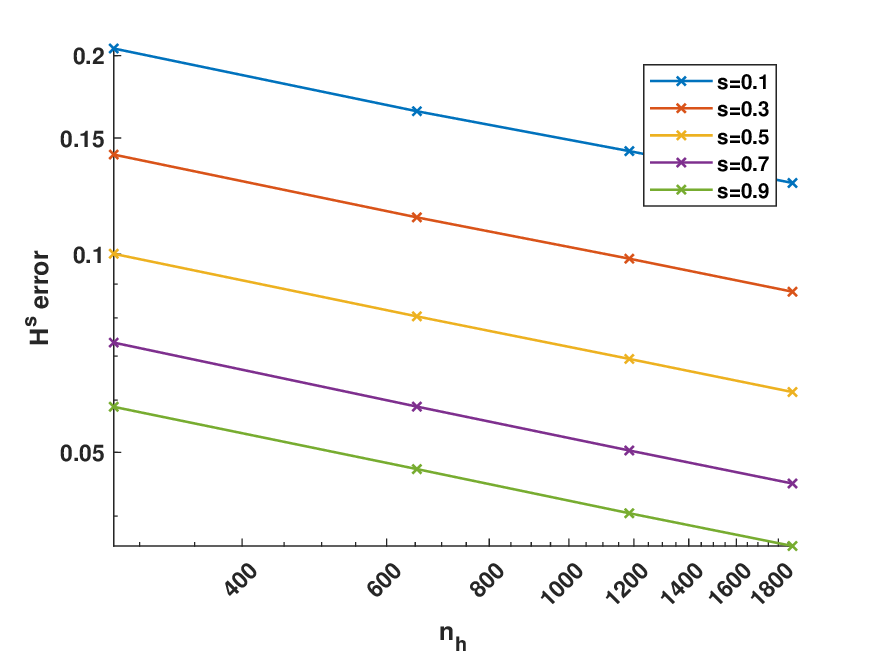}
    \end{minipage}
    \hfill
    \begin{minipage}{0.49\textwidth}
    \centering
    \includegraphics[width=\textwidth]{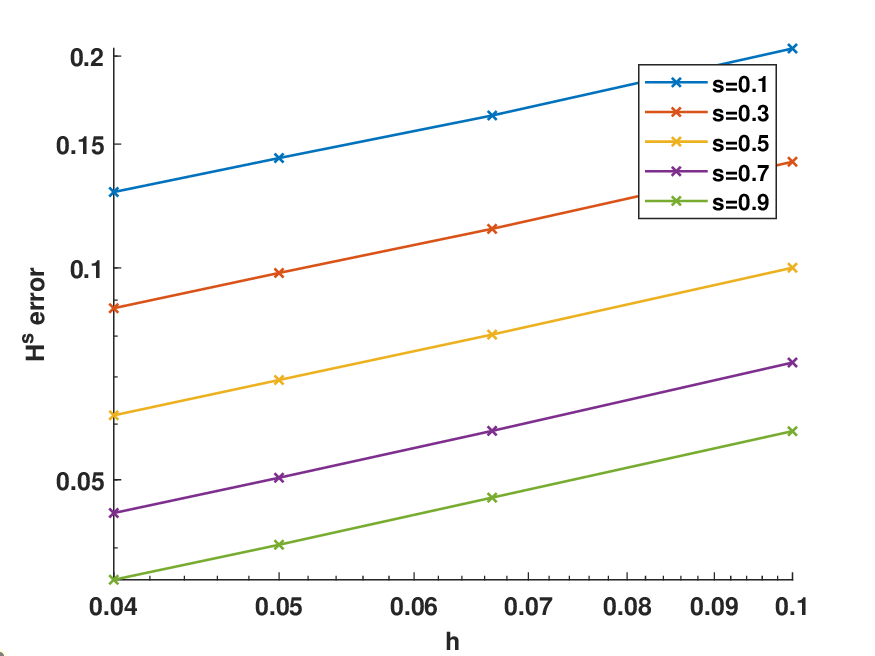}
    \end{minipage}
    \caption{Convergence results for $d=2$ corresponding to Table \ref{tab:errores_2d} in logarithmic scale.} \label{fig:conv-rates-d=2}
\end{figure}

\begin{figure}[!htbp]
    \centering
    \begin{minipage}{0.49\textwidth}
    \centering
    \includegraphics[width=\textwidth]{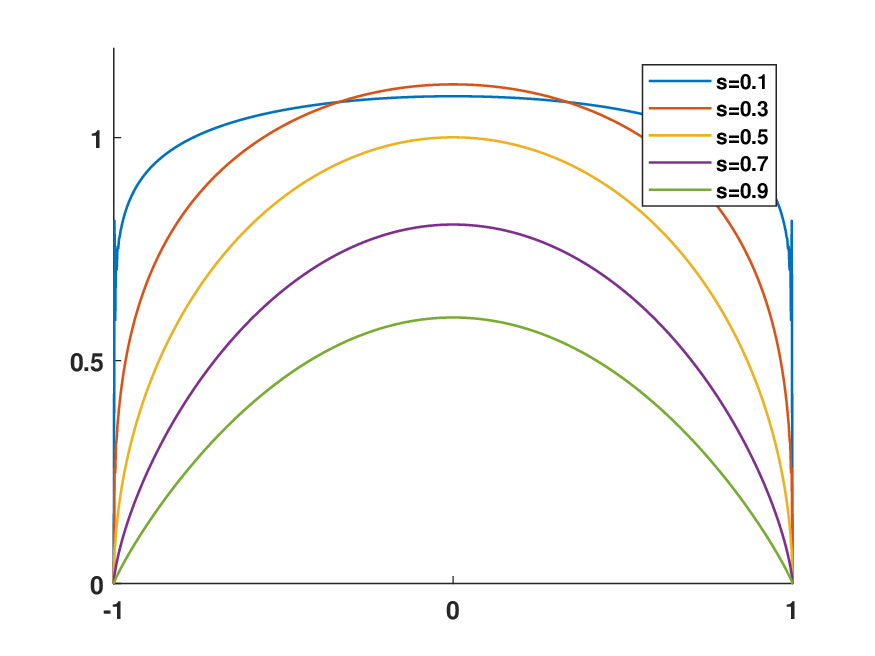}
    
    \end{minipage}
    \hfill
    \begin{minipage}{0.49\textwidth}
    \centering
    \includegraphics[width=\textwidth]{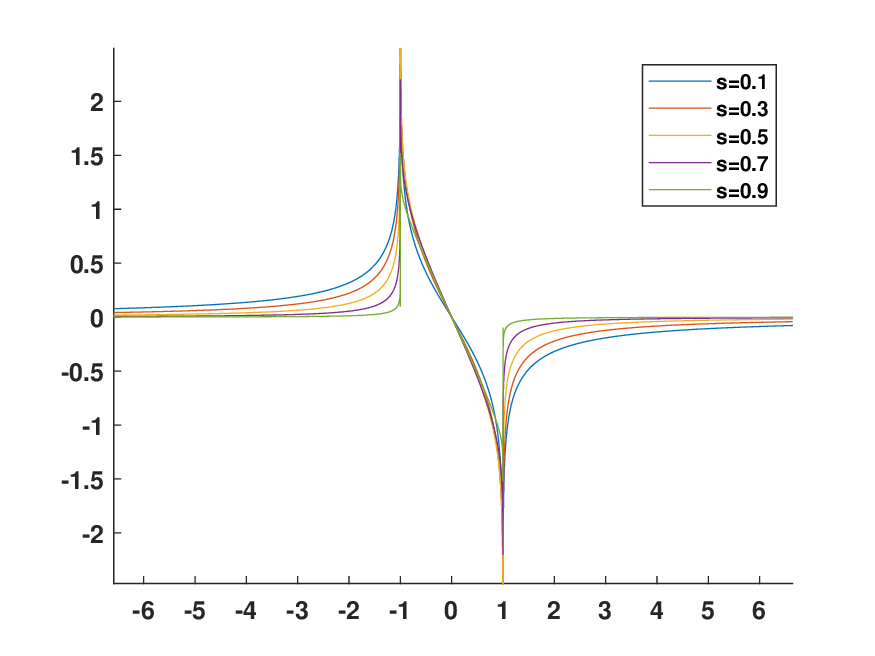}
    \end{minipage}
    \caption{Here $f\equiv 1$ and $\Omega=(-1,1)$. The left panel shows the computed pressures for different values of $s$. The right panel displays the computed fluxes for different values of $s$. The mesh size parameter is $h=0.022$.} \label{fig:soluciones_1d}
\end{figure}

\begin{figure}[!htbp]
    \centering
    \begin{minipage}{0.32\textwidth}
        \includegraphics[width=\textwidth]{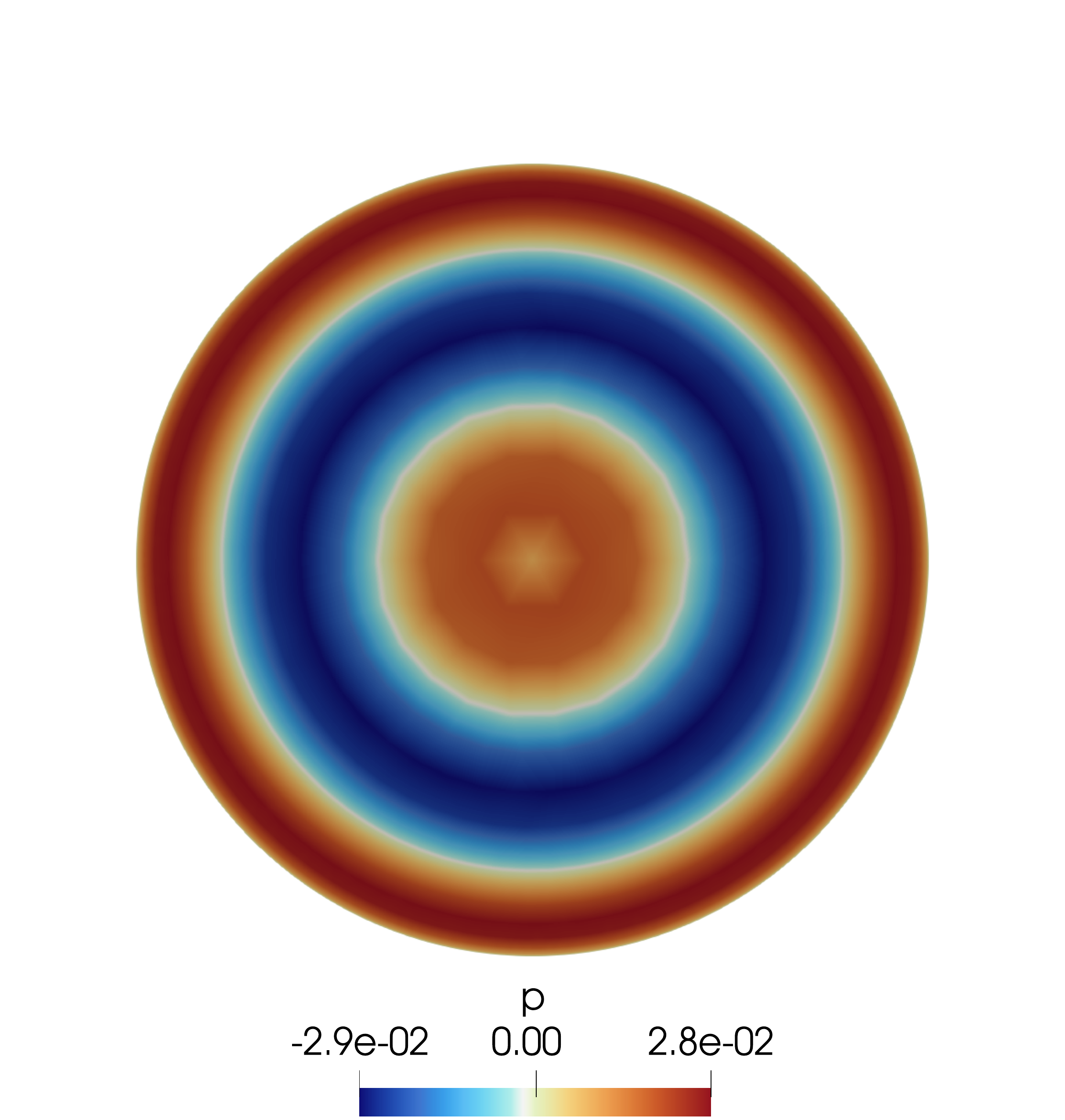}
    \end{minipage}
    \hfill
    \begin{minipage}{0.32\textwidth}
        \includegraphics[width=\textwidth]{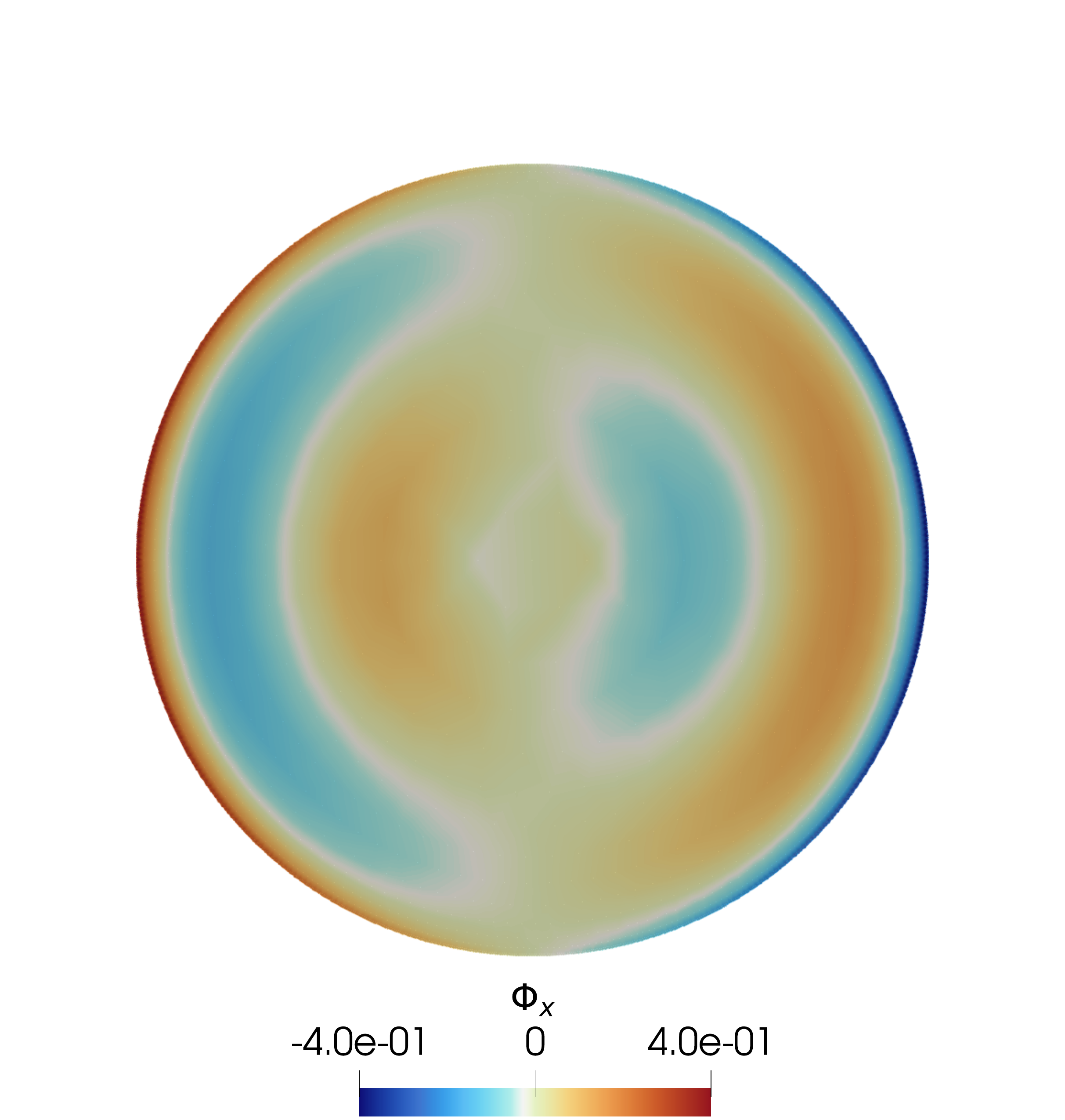}
    \end{minipage}
    \hfill
        \begin{minipage}{0.32\textwidth}
        \includegraphics[width=\textwidth]{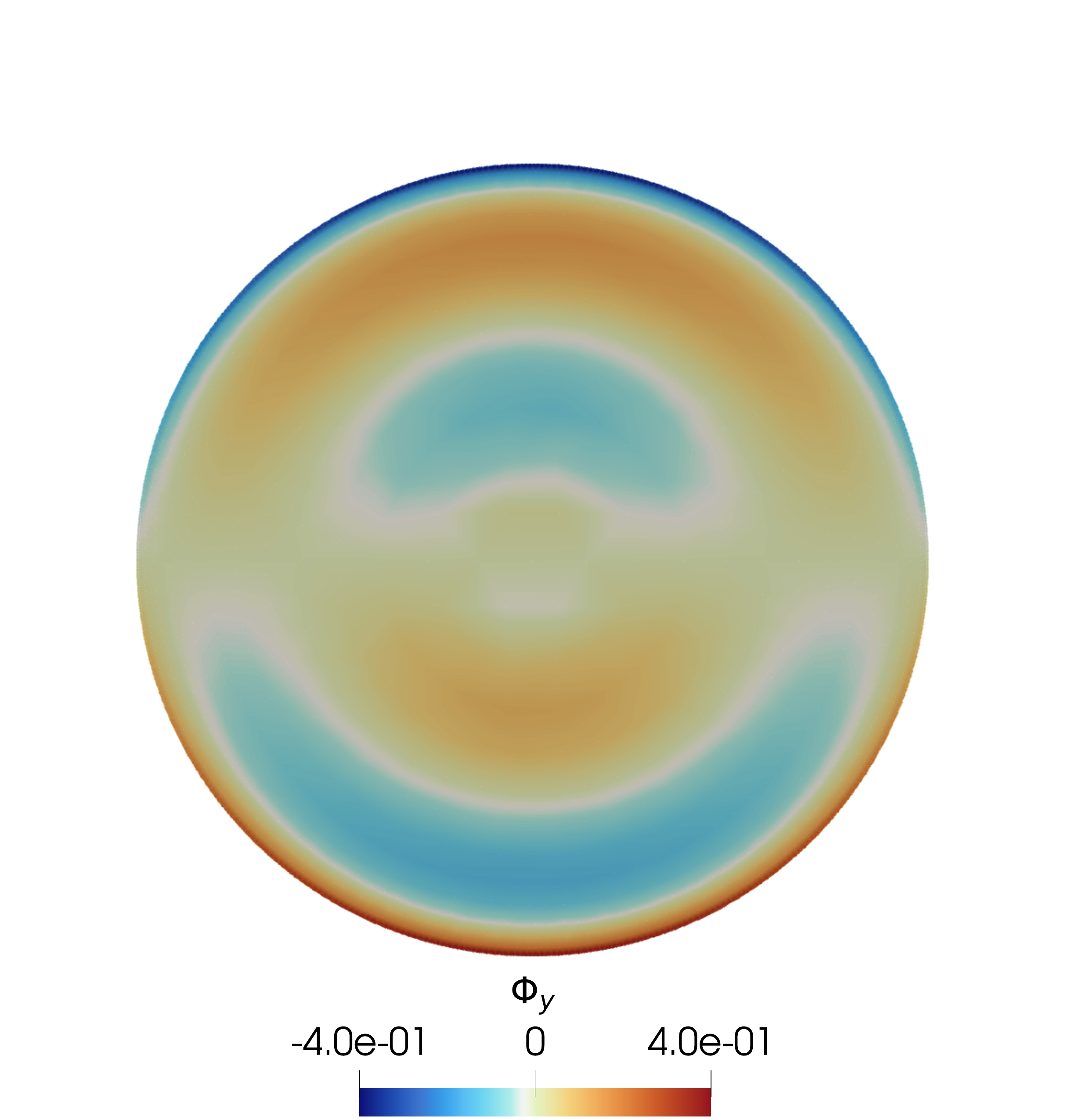}
    \end{minipage}
    \caption{Here, $f(x)= P_{s,0,3}(2|x|^2-1)$ as in \eqref{eq:Jacobi-pols}, $s=0.5$, $\Omega=B(0,1)$, and $B_H=B(0,2)$ as mentioned before. The left panel shows the computed pressure and the center and right ones display the components of the computed flux inside $\Omega$. We highlight the singular behavior of the flux near the boundary of $\Omega$. The mesh contains 11406  elements with 2401 degrees of freedom for the pressure and $2\times5822$ degrees of freedom for the flux.}     \label{fig:soluciones_2d} 
\end{figure}

For the pressure $p$, the computed orders of convergence for different values of $s$ in one and two dimensions are respectively presented in Tables \ref{tab:errores_1d} and \ref{tab:errores_2d} and Figure \ref{fig:conv-rates-d=2}. Furthermore, Figure \ref{fig:soluciones_1d} displays the computed pressures and fluxes for different values of $s$ in one dimension, and Figure \ref{fig:soluciones_2d} shows an example of computed solutions in two dimensions with $s=0.5$ and $n=3$ in Theorem \ref{thm:explicit_solutions}. The outcomes of our computational experiments are consistent with Propositions \ref{cor:order-convergence} and \ref{prop:Aubin-Nitsche}. 

\begin{table}[!htbp]
    \centering
    \begin{tabular}{c|c|c}
        Value of $s$ & $H^s$-order & $L^2$-order  \\
        \hline
        0.1 & 0.4691 & 0.5949\\
        0.2 & 0.4956  & 0.6444\\
        0.3 & 0.5000  & 0.7968\\
        0.4 & 0.5004  & 0.9236\\
        0.5 & 0.5005  & 1.0012\\
        0.6 & 0.5005  & 0.9966\\
        0.7 & 0.5009 & 0.9928\\
        0.8 & 0.5014  & 0.9952\\
        0.9 & 0.5017  & 1.0045\\
    \end{tabular}
    \caption{Order of convergence for the pressure $p$ in problem \eqref{eq:Darcy-torsion} in the one dimensional case.}
    \label{tab:errores_1d}
\end{table}
\begin{table}[!htbp]
    \centering
    \begin{tabular}{c|c|c|c|c}
        Value of $s$ & $H^s$-order in $h$ & $L^2$-order in $h$ & $H^s$-order in $n_h$ & $L^2$-order in $n_h$ \\
        \hline
        0.1 & 0.4985 & 0.5869 & -0.2430 & -0.2861 \\
        0.2 & 0.4959  & 0.6817 & -0.2417 & -0.3323\\
        0.3 & 0.5170 & 0.8309 & -0.2520 & -0.4050\\
        0.4 & 0.5314  & 0.9208 & -0.2590 & -0.4488\\
        0.5 & 0.5187  & 0.9989 & -0.2528 & -0.4869\\
        0.6 & 0.5189  & 1.1164 & -0.2529 & -0.5442\\
        0.7 & 0.5175 & 1.2247 & -0.2523 & -0.5969\\
        0.8 & 0.5127  & 1.1923 & -0.2499 & -0.5811\\
        0.9 & 0.5131  & 1.0946 & -0.2501 & -0.5336\\
    \end{tabular}
    \caption{Order of convergence for the pressure $p$ in problem \eqref{eq:Darcy-torsion} the two dimensional case.}
    \label{tab:errores_2d}
\end{table}

\subsection{Dependence on $H$}
Our second set of experiments deals with the dependence of discrete solutions to \eqref{eq:Darcy-torsion} with respect to the computational domain diameter $H$.
For a fixed $h$ and different values of $s$, we compute errors for different values of $H$. We recall that our approach extends the discrete flux $\vPhi_h$ by zero outside $B_H$, so that $H$ effectively acts as a truncation parameter. Consequently, an improvement of the error is expected as $H$ increases. Naturally, this will increase the number of elements of the mesh and, therefore, the CPU time used in the computations. Nevertheless, Proposition \ref{cor:order-convergence} suggests that significant improvements in the error are not be expected once 
$H$ is sufficiently large.
Our results, displayed in Table \ref{tab:exp2-varR}, are in good agreement with this heuristic idea.

\begin{table}[!htbp]
\label{tab:exp2-varR}
    \centering
    \begin{minipage}{0.48\textwidth}
        \centering
        \begin{tabular}{c|c|c|c}
            \textbf{$H^s$-error} & $s = 0.2$ & $s=0.5$ & $s=0.8$ \\
            \hline
            $H=0.40$ & 0.1444 & 0.0724 & 0.0445 \\
            $H=0.70$ & 0.1299 & 0.0712 & 0.0444 \\
            $H=1.50$ & 0.1219 & 0.0700 & 0.0435 \\
            $H=1.70$ & 0.1215 & 0.0699 & 0.0435 \\
            $H=1.90$ & 0.1212 & 0.0701 & 0.0436 \\
            $H=2.02$ & 0.1207 & 0.0696 & 0.0431 \\
        \end{tabular}
    \end{minipage}
    \caption{$H^s$-error for the pressure in \eqref{eq:Darcy-torsion} in $d=2$ over quasi-uniform meshes with  size $h=0.07$ and different computational domain diameters.}
\end{table}

\subsection{Graded meshes}
As a final example, we test the convergence rates of the method when using graded meshes. Theorem \ref{thm:regularidad_con_pesos} suggests that improved convergence rates could be achieved by using a priori adapted meshes in the spirit of \cite{acosta2017fractional}. Given $H>0$,
we consider a mesh $\Th(B_H)$ satisfying the following requirement for elements contained in $\Omega$:
\begin{equation}
 \label{eq:graded_mesh}\tag{G}
 h_T \approx 
 \begin{cases}
h\,d(T,\partial\Omega)^{\frac{1}{2}} & \mbox{if } T\cap\partial\Omega=\emptyset, \\     h^2 & \mbox{if } T\cap\partial\Omega \neq \emptyset.
 \end{cases}
\end{equation}
Over $\Omega= B(0,1)$ such a mesh can be obtained in the following way. Consider an integer $N>0$ and the sequence $r_j=(1-\frac{j}{N})^2$ for $1\leq j\leq N$. Let $\Th(B_H)$ be the union of all uniform meshes with mesh size $h_j=r_{j}-r_{j-1}$ in the domains $\{x\in B(0,1):r_{j-1}<|x|< r_j\}\subset\Omega$ for $1\leq j\leq N$. This construction ensures that conditions \eqref{eq:graded_mesh} are satisfied with $h \simeq 1/N$. 

The improved orders of approximation for different values of $s$ are displayed in Table \ref{tab:orden_convergencia_graduadas} for the torsion problem \eqref{eq:Darcy-torsion}.
As expected from the weighted regularity in Theorem \ref{thm:regularidad_con_pesos}, first-order convergence for the pressure is attained over these meshes.

\begin{table}[!htbp]
    \centering
    \begin{tabular}{c|c}
        Value of s & $H^s$ order in $h$ \\
        \hline
        0.1 & 0.9601  \\
        0.2 & 0.9873  \\
        0.3 & 1.0059 \\
        0.4 & 1.0181  \\
        0.5 & 1.0269  \\
        0.6 & 1.0349  \\
        0.7 & 1.0437 \\
        0.8 & 1.0553  \\
        0.9 & 1.0702  \\
    \end{tabular}
    \caption{Order of convergence for the pressure $p$ using graded meshes. Here, $f\equiv1$ and $\Omega=B(0,1)$.}
    \label{tab:orden_convergencia_graduadas}
\end{table}

\section*{Acknowledgments}
The authors have been supported in part by Fondo Clemente Estable grant 172393 and program MATH-AmSud 23MATH-06. 

\bibliographystyle{abbrv}
\bibliography{darcy}
\end{document}